\documentclass[11pt]{amsart}

\usepackage[hmargin=1.0in,height=8.3in]{geometry}
\usepackage{amsthm,amsmath,epsfig,subfigure,amsfonts,amssymb,bm,xcolor,amsbsy,times}
\usepackage{bm}

\usepackage{ifpdf}
\usepackage{color}
\definecolor{webgreen}{rgb}{0,.5,0}
\definecolor{webbrown}{rgb}{.8,0,0}
\definecolor{emphcolor}{rgb}{0.95,0.95,0.95}

\usepackage{hyperref}
\hypersetup{%
          colorlinks=true,
          linkcolor=blue,
          filecolor=webbrown,
          citecolor=webbrown,
          breaklinks=true}
\ifpdf \hypersetup{pdftex,
            pdfstartview=FitH, 
            bookmarksopen=true,
            bookmarksnumbered=true
} \else \hypersetup{dvips} \fi
\usepackage[square,sort,comma,numbers]{natbib}

\numberwithin{equation}{section}

\newtheorem{theorem}{Theorem}[section]
\newtheorem{proposition}{Proposition}[section]

\newtheorem{remark}{Remark}[section]
\newtheorem{lemma}{Lemma}[section]

\newtheorem{Assumption}{Assumption}[section]
\newtheorem{definition}{Definition}[section]
\numberwithin{remark}{section} 
\numberwithin{proposition}{section}
\numberwithin{corollary}{section}

\allowdisplaybreaks[1]

\title[ ]{Optimal Entry and Consumption Under Habit Formation}

\author[  ]{Yue Yang}
\address[Y. Yang]{Department of Applied Mathematics, The Hong Kong Polytechnic University, Hung Hom, Kowloon, Hong Kong.}
\email{yue.yy.yang@connect.polyu.hk}
\author[  ]{Xiang Yu}
\address[X. Yu]{Department of Applied Mathematics, The Hong Kong Polytechnic University, Hung Hom, Kowloon, Hong Kong.}
\email{xiang.yu@polyu.edu.hk}

\begin{document}

\begin{abstract}
This paper studies a composite problem involving the decision making of the optimal entry time and dynamic consumption afterwards. In stage-1, the investor has access to full market information subjecting to some information costs and needs to choose an optimal stopping time to initiate stage-2; in stage-2, the investor terminates the costly full information acquisition and starts dynamic investment and consumption under partial observations of free public stock prices. The habit formation preference is employed, in which the past consumption affects the investor's current decisions. By using the stochastic Perron's method, the value function of the composite problem is proved to be the unique viscosity solution of some variational inequalities.
\end{abstract}

\maketitle

\section{Introduction}
We consider a simple model to incorporate information costs in a continuous time portfolio-consumption problem. In particular, we study a two-stage composite problem under complete and incomplete filtrations sequentially. The drift process of the stock price is assumed to be of the Ornstein-Uhlenbeck type. In the first stage from the initial time, the investor needs to pay information costs to access the full market information generated by both drift and stock price processes to update their dynamic distributions and decide the optimal time to enter the second stage. The information costs may refer to search cost, storage cost, communication cost, investor's attention cost or other service costs. We consider the simple linear information cost in the present paper, which is modeled by a constant cost rate and will be subtracted directly from the investor's wealth amount. That is, the longer the first stage is, the higher information costs the investor needs to afford. Some previous work have addressed impacts of information costs to optimal investment from different perspectives, see \cite{kang1997there}, \cite{portes2005determinants}, \cite{ahearne2004information} and \cite{keppo2018investment}. In our first stage, the mathematical problem becomes an optimal stopping problem under the complete market information filtration. The second stage starts from the chosen entry time and the investor terminates the full observations of the drift process. Instead, the investor starts to dynamically choose the investment and consumption based on the prior data inputs and the free partial observations of the stock price, which can be formulated as an optimal control problem under incomplete information filtration. As the value function of the interior control problem depends on the stopping time and data inputs of the drift process, the exterior problem can be interpreted as to wait in an optimal way so that the input values can achieve the maximum of the interior functional.

Portfolio optimization under partial observations have been extensively studied in past decades, see a few examples in \cite{lakner1998optimal,xia2001learning,brendle2006portfolio,monoyios2009optimal,bjork2010optimal,BoLY} with different financial motivations. As illustrated in these work, the value function under incomplete information filtration is strictly lower than the counterpart under full information filtration and this gap is usually regarded as the value of information. The present paper attempts to study partial observations from a different perspective that the full market information is available but costly because more data, services and personal attentions are involved. The information costs may change the investor's attitude towards the usage of full observations because it is no longer true that the more information he observes, the higher profit he can attain. Moreover, from some previous work on partial observations, we know that the value function eventually depends on the given initial input of the random factor such as the drift process. As in \cite{lakner1998optimal,brendle2006portfolio}, it is conventionally assumed that the initial data of the unobservable drift is a Gaussian random variable so that the Kalman-Bucy filtering can be applied. We take this input into account and consider a model that the investor can wait and dynamically update the distribution of inputs using the full market information subjecting to information costs. We can show that starting sharp from the initial time to invest and consume under incomplete information may not be optimal. 

On the other hand, the habit formation has become a new paradigm for modelling preferences on consumption rate in recent years, which can better match with some empirical observations, see \cite{constantinides1990habit,mehra1985equity}. The literature suggests that the past consumption pattern may enforce a continuing impact on individual's current consumption decisions. In particular, the linear habit formation preference has been widely accepted, in which there exists an index term that stands for the accumulated consumption history. The habit formation preference has been well studied by \cite{detemple1992optimal,englezos2009utility,munk2008portfolio} in complete market models and by \cite{yu2015utility,yu2014optimal} in incomplete market models. It is noted that the utility function is decreasing in the habit level. In the present paper, we assume that there is no consumption during stage-1 and the investor starts to gain consumption habit only in stage-2. Therefore, an early entry time to stage-2 may not be the optimal decision because the investor has longer time to develop a much higher habit level. This is our second motivation to investigate the exterior optimal entry time problem to see whether longer waiting can benefit the investor more as the resulting habit level can be much lower that may lead to a higher value function.

We show that the value function of the composite problem is the unique viscosity solution to some variational inequalities. To this end, we can choose to apply either classical Perron's method or the stochastic version of Perron's method introduced in \cite{bayraktar2011probabilistic}. For classical Perron's method, to establish the equivalence between the value function and the viscosity solution, we have to either prove the dynamic programming principle or upgrade the global regularity of the solution and prove the verification theorem. The convexity (concavity) of the value function with respect to the state variable is usually crucial in some standard arguments to conclude the global regularity. However, this property is not clear in our composite problem, see Remark \ref{convexconcave} for details. The global regularity of the value function is not guaranteed, and the direct verification proof for our exterior problem becomes difficult. Therefore, we choose the stochastic Perron's method instead, which allows us to show the equivalence between the value function and the viscosity solution without global regularity. For some related works on optimal stopping using viscosity solution, we refer to \cite{reikvam1998viscosity} and \cite{pham1997optimal}. See also some recent work on stochastic control problems using stochastic Perron's method among \cite{bayraktar2011probabilistic,bayraktar2013stochastic,bayraktar2014stochastic,bayraktar2015stochastic,sirbu2014stochastic,LeeYZ}. One important step to complete the argument of stochastic Perron's method is the comparison principle of the associated variational inequalities, which is established in the present paper.

The rest of the paper is organized as follows: Section \ref{sec-2} introduces the market model and the habit formation preference and formulates the 2-stage optimization problem. Section \ref{sec-3} gives the main result of the interior utility maximization problem with habit formation and partial observations. Section \ref{sec-4} studies the exterior optimal entry problem with linear information costs. Using the stochastic Perron's method, we show that the value function of the composite problem is the unique viscosity solution of some variational inequalities. Some auxiliary results and proofs are reported in Appendix \ref{appendix5.1} and \ref{appendix5.2}.

\section{Mathematical Model and Preliminaries}\label{sec-2}
\subsection{Market model}

Given the probability space $(\Omega,\mathbb{F},\mathbb{P})$ with full information filtration $\mathbb{F}=(\mathcal{F}_{t})_{0\leq t\leq T}$ that satisfies the usual conditions, we consider the market with one risk-free bond and one risky asset over a finite time horizon $[0,T]$. It is assumed that the bond process satisfies $S^{0}_{t}\equiv1$, for $t\in[0,T]$, which amounts to the standard change of num\'eraire.

The stock price $S_{t}$ satisfies
\begin{align}\label{dst}
      dS_{t}=\mu_{t}S_{t}dt+\sigma_{S}S_{t}dW_{t}, \ \ \  0\leq t\leq T,
\end{align}
with $S_{0}=s>0$. Some empirical studies such as \cite{brennan2010persistence,campbell1997econometrics,fama1989business,poterba1988mean} have observed that the drift process of many risky assets follows the so-called mean reverting diffusion. We also consider here that the drift process $\mu_t$ in \eqref{dst} satisfies the Ornstein-Uhlenbeck SDE by
\begin{align}\label{realmu}
      d\mu_{t}=-\lambda(\mu_{t}-\bar{\mu})dt+\sigma_{\mu}dB_{t}, \ \ \  0\leq t\leq T.
\end{align}
Here, $(W_{t})_{0\leq t\leq T}$ and $(B_{t})_{0\leq t\leq T}$ are $\mathcal{F}_{t}$-adapted Brownian motions with correlation coefficient $\rho\in[-1,1]$. For simplicity, the initial value $\mu_0$ of the drift is a given constant. We assume that market coefficients $\sigma_{S}$, $\lambda$, $\bar{\mu}$ and $\sigma_{\mu}$ are given nonnegative constants based on calibrations from historical data.

It is assumed that the investor starts with initial wealth $x(0)=x_0>0$ at time $t=0$. Also, starting from the initial time $t=0$, the access to the full market information $\mathcal{F}_t$ generated by $W$ and $B$ incurs information costs $\kappa t$, where $\kappa>0$ is the constant cost rate per unit time. The information costs may refer to storage cost, search cost, communication cost, investor's attention cost or other service costs to fully observe the market information $\mathcal{F}_t$. Moreover, to simplify the mathematical problem, it is assumed that starting from $t=0$ to a chosen stopping time $\tau$, the investor purely waits and updates dynamic distributions of processes $\mu_t$ and $S_t$ and does not invest and consume at all. This assumption makes sense as long as the value of the optimal entry time $\tau$ is short in the model. The dynamic wealth process after the information costs at time $t$ is simply given by a deterministic function $x(t)=x_0-\kappa t$ for any $t\leq \tau$.

As the full market information filtration is costly, the investor needs to choose optimally choose an $\mathcal{F}_{t}$-adapted stopping time $\tau$ to terminate the full information acquisition and enter the second stage. From the chosen stopping time $\tau$, he switches to the partial observations filtration $\mathcal{F}_t^S=\mathcal{F}_{\tau}\bigvee\sigma(S_u:\tau\leq u\leq t)$ for $\tau\leq t\leq  T$, which is the union of the sigma algebra $\mathcal{F}_{\tau}$ and the natural filtration generated by the stock price $S$ up to time $t$. Moreover, for any time $\tau\leq t\leq T$, the investor chooses a dynamic consumption rate $c_{t}\geq 0$ and decides the amounts $\pi_{t}$ of his wealth to invest in the risky asset and the rest in the bond. Without paying information costs, the drift process $\mu_{t}$ and Brownian motions  $W_{t}$ and $B_{t}$ are no longer observable for $t\geq \tau$. Therefore, the investment-consumption pair $(\pi_t, c_t)$ is only assumed to be adapted to the partial observation filtration $\mathcal{F}_t^S$ for $\tau\leq t\leq T$. Recall that at the entry time $\tau$, the investor only has wealth $x(\tau)=x_0-\kappa \tau$ left. Under the incomplete filtration $\mathcal{F}_t^S$, the investor's total wealth process $\hat{X}_{t}$ can be written as
    \begin{align}\label{dxt}
      d\hat{X}_{t}=(\pi_{t}\mu_{t}-c_{t})dt+\sigma_{S}\pi_{t}dW_{t}, \ \ \  \tau\leq t\leq T,
    \end{align}
with the initial value $\hat{X}_{\tau}=x(\tau)=x_{0}-\kappa \tau>0$. Note that $W_t$ is no longer a Brownian motion under the partial observations filtration $\mathcal{F}_t^S$, we have to apply the Kalman-Bucy filtering and consider the \emph{Innovation Process} defined by
    \begin{align*}
      d\hat{W}_{t}:=\frac{1}{\sigma_{S}}\Big[(\mu_{t}-\hat{\mu}_{t})dt+\sigma_{S}dW_{t}\Big]=\frac{1}{\sigma_{S}}\Big(\frac{dS_{t}}{S_{t}}-\hat{\mu}_{t}dt\Big), \ \ \  \tau\leq t\leq T,
    \end{align*}
    which is a Brownian motion under $\mathcal{F}^{S}_{t}$. The best estimation of the unobservable drift process $\mu_t$ under $\mathcal{F}_t^S$ is the conditional expectation process $\hat{\mu}_{t}=\mathbb{E}\Big[\mu_{t}\Big{|}\mathcal{F}^{S}_{t}\Big]$, for $\tau\leq t\leq T$ with the initial input $\hat{\mu}_{\tau}=\mu_{\tau}$, $\mathbb{P}$-a.s., at the stopping time $\tau$ where the distribution of $\mu_{\tau}$ is determined via \eqref{realmu} by paying information costs up to $\tau$. By standard Kalman-Bucy filtering (see equation (18) of \cite{brendle2006portfolio} or equation (21) of \cite{monoyios2009optimal}), $\hat{\mu}_t$ satisfies the SDE
    \begin{align}\label{dhatmut}
      d\hat{\mu}_{t}=-\lambda(\hat{\mu}_{t}-\bar{\mu})dt+\left(\frac{\hat{\Sigma}(t)+\sigma_{S}\sigma_{\mu}\rho}{\sigma_{S}}\right)d\hat{W}_{t}, \ \ \  \tau\leq t\leq T,
    \end{align}
    with $\hat{\mu}_{\tau}=\mu_{\tau}$, $\mathbb{P}$-a.s.. Moreover, the conditional variance $\hat{\Sigma}(t)=\mathbb{E}\Big[(\mu_{t}-\hat{\mu}_{t})^2\Big{|}\mathcal{F}^{S}_{t}\Big]$ satisfies the deterministic Riccati ODE (see equation (19) of \cite{brendle2006portfolio} or equation (23) of \cite{monoyios2009optimal})
    \begin{align}\label{domega}
      \frac{d\hat{\Sigma}(t)}{dt}=-\frac{1}{\sigma_{S}^2}\hat{\Sigma}^2(t)+\Big(-\frac{2\sigma_{\mu}\rho}{\sigma_{S}}-2\lambda\Big)\hat{\Sigma}(t)+(1-\rho^2)\sigma_{\mu}^2, \ \ \  \tau\leq t\leq T,
    \end{align}
    with the initial value $\hat{\Sigma}(\tau)= \mathbb{E}\Big[(\mu_{\tau}-\hat{\mu}_{\tau})^2\Big{|}\mathcal{F}^{S}_{t}\Big] = 0$ in view of $\hat{\mu}_{\tau} = \mu_{\tau}$, $\mathbb{P}$-a.s.. It can be solved explicitly as
    \begin{align*}
      \hat{\Sigma}(t)=\sqrt{k}\sigma_{S}\frac{k_{1}\exp(2(\frac{\sqrt{k}}{\sigma_{S}})t)+k_{2}}{k_{1}\exp(2(\frac{\sqrt{k}}{\sigma_{S}})t)-k_{2}}-\left(\lambda+\frac{\sigma_{\mu}\rho}{\sigma_{S}}\right)\sigma_{S}^2, \ \ \  \tau\leq t\leq T,
    \end{align*}
    where $k=\lambda^2\sigma_{S}^2+2\sigma_{S}\sigma_{\mu}\lambda\rho+\sigma_{\mu}^2$, $k_{1}=\sqrt{k}\sigma_{S}+(\lambda\sigma_{S}^2+\sigma_{S}\sigma_{\mu}\rho)$ and $k_{2}=-\sqrt{k}\sigma_{S}+(\lambda\sigma_{S}^2+\sigma_{S}\sigma_{\mu}\rho)$.

    For the second stage dynamic control problem, we employ the habit formation preference. In particular, we denote $Z_{t}:= Z(c_{t})$ as \textit{habit formation} process or \textit{the standard of living} process, which describes the consumption habits level. It is assumed conventionally that the accumulative reference $Z_{t}$ satisfies the recursive equation (see \cite{detemple1992optimal}) that $dZ_{t}=(\delta(t)c_{t}-\alpha(t)Z_{t})dt$, $\tau\leq t\leq T$, where $Z_{\tau}=z_{0}\geq 0$ is called the \emph{initial consumption habit} of the investor. Equivalently, we have
    \begin{align*}
      Z_{t}=z_{0}e^{-\int^{t}_{\tau}\alpha(u)du}+\int^{t}_{\tau}\delta(u)e^{-\int^{t}_{u}\alpha(s)ds}c_{u}du, \ \ \  \tau\leq t\leq T,
    \end{align*}
    which is the exponentially weighted average of the initial habit and the past consumption. Here, the deterministic discount factors $\alpha(t)\geq 0$ and $\delta(t)\geq 0$ measure, respectively, the persistence of the past level and the intensity of consumption history. We are interested in \emph{addictive habits} in the present paper, namely it is required that the investor's current consumption strategies shall never fall below the level of standard of living that $c_{t}\geq Z_{t}$ a.s., for $\tau\leq t\leq T$.

    Under the partial observation filtration $(\mathcal{F}^{S}_{t})_{\tau\leq t\leq T}$, the stock price dynamics (\ref{dst}) can be rewritten by $dS_{t}=\hat{\mu}_{t}S_{t}dt+\sigma_{S}S_{t}d\hat{W}_{t}$ and the wealth dynamics (\ref{dxt}) can be rewritten as $d\hat{X}_{t}=(\pi_{t}\hat{\mu}_{t}-c_{t})dt+\sigma_{S}\pi_{t}d\hat{W}_{t}$, $\tau\leq t\leq T$. To facilitate the formulation of the stochastic control problem and the derivation of the dynamic programming equation, for any $t\in [0,T]$, we denote $\mathcal{A}_{t}(x)$ the time-modulated admissible set of the pair of investment and consumption process $(\pi_{s},c_{s})_{t\leq s\leq T}$ with the initial wealth $\hat{X}_t=x$, which is $\mathcal{F}^{S}_{s}$-progressively measurable and satisfies the integrability conditions  $\int^{T}_{t}\pi_{s}^2ds<+\infty$, $\mathbb{P}$-a.s. and $\int^{T}_{t}c_{s}ds<+\infty$, $\mathbb{P}$-a.s..
    with the addictive habit formation constraint that $c_{s}\geq Z_{s}$, $\mathbb{P}$-a.s., $t\leq s\leq T$. Moreover, no bankruptcy is allowed, i.e., the investor's wealth remains nonnegative, i.e. $\hat{X}_{s}\geq 0$, $\mathbb{P}$-a.s., $t\leq s\leq T$.

\subsection{Problem formulation}

The two-stage optimal decision making problem is formulated as the composite problem involving the optimal stopping and the stochastic control afterwards, which is defined by
   \begin{align}\label{0timev}
   \begin{aligned}
   \widetilde{V}(0,\mu_0;x_0, z_0) &:=\sup_{\tau\geq 0}\mathbb{E}\left[
   \underset{(\pi,c)\in\mathcal{A}_{\tau}(x_0-\kappa\tau)}{\textrm{esssup}}\mathbb{E}\left[\int^T_\tau \frac{(c_{s}-Z_{s})^p}{p}ds\Big{|}\mathcal{F}_\tau^{S}\right]\right].
   \end{aligned}
   \end{align}
In particular, starting from the chosen stopping time $\tau$, we are interested in the utility maximization on consumption with habit formation, in which the power utility function $U(x)=x^p/p$ is defined on the difference $c_t-Z_t$. To simplify the presentation, we only consider in the present paper that the risk aversion coefficient $p<0$. The indirect utility process of the interior control problem is denoted by
    \begin{align*}
    &\widehat{V}(t,x_0-\kappa t, z_0, \mu_t; 0) :=
    \underset{(\pi,c)\in\mathcal {A}_{t}(x_0-\kappa t)}{\textrm{esssup}}\mathbb{E}\left[\int^T_t \frac{(c_{s}-Z_{s})^p}{p}ds\Big{|}\mathcal{F}_t^{S}\right]\\
    =&
    \underset{(\pi,c)\in\mathcal {A}_{t}(x_0-\kappa t)}{\textrm{esssup}}\mathbb{E}\left[\int^T_t \frac{(c_{s}-Z_{s})^p}{p}ds\Big{|}\hat{X}_{t}=x_0-\kappa t,\hat{\mu}_{t}=\mu_{t},Z_{t}=z_0;\hat{\Sigma}(t)=0\right].
    \end{align*}

To determine the exterior optimal stopping time, we need to maximize over the inputs of values $\tau$, $\hat{X}_{\tau}$, $Z_{\tau}$ and $\hat{\mu}_{\tau}$. Recall that the investor does not manage his investment and consumption before $\tau$, it follows that $\hat{X}_{\tau}=x_0-\kappa\tau$, $Z_{\tau}=z_0$ and $\hat{\Sigma}(\tau)=0$ can all be taken as parameters instead of variables. That is, $\mu_{\tau}=\hat{\mu}_{\tau}$ is the only random input and we can regard $\mu_t$ as the only underlying state process. Therefore, the dynamic counterpart of \eqref{0timev} is defined by
   \begin{align}\label{tilde2v}
     \widetilde{V}(t,\eta;x_0-\kappa t,z_0): = \underset{\tau\geq t}{\textrm{esssup}} \ \mathbb{E}\left[\underset{(\pi,c)\in\mathcal{A}_{\tau}(x_0-\kappa\tau)}{\textrm{esssup}}\ \mathbb{E}\left[\int^T_\tau \frac{(c_{s}-Z_{s})^p}{p})ds\Big{|}\mathcal{F}_\tau^{S}\right]\Bigg{|}\mu_t=\eta\right].
   \end{align}

\begin{remark}
We focus on the case $p<0$ in the present paper because functions $A(t,s)$, $B(t,s)$ and $C(t,s)$ as solutions to some future ODEs (\ref{atts}), (\ref{btts}), (\ref{ctts}) are all bounded and the utility $U(x)$ is also bounded from above, which can significantly simplify the proof of the verification result in Theorem \ref{thm} and the proof of comparison results in Proposition \ref{compp}. The other case $0<p<1$ can essentially be handled in a similar way. However, as the process $\hat{\mu}_t$ in \eqref{dhatmut} is unbounded and functions $A(t,s)$, $B(t,s)$ and $C(t,s)$ may explode at some $t\in [0,T]$, one needs some additional parameter assumptions to guarantee integrability conditions and martingale properties in the proofs of some main results.
\end{remark}

\begin{Assumption}\label{f}
According to Remark \ref{effdomain} for the interior control problem, it is assumed from this point onwards that $x_0-\kappa t>z_0 m(t)$ for any $0\leq t\leq T$, i.e. the initial wealth is sufficiently large after paying information costs such that the interior control problem is well defined for any $0\leq t\leq T$, where $m(t)$ is defined by
    \begin{align}\label{mt}
      m(t)=\int^T_t\exp\left(\int^s_t(\delta(v)-\alpha(v))dv\right)ds, \ \ 0\leq t\leq T.
    \end{align}
    We note that $m(t)$ in \eqref{mt} represents the cost of subsistence consumption per unit of standard of living at time $t$ because the interior control problem is solvable if and only if $\hat{X}_t^{*}\geq m(t)Z_t$, $0\leq t\leq T$, see Lemma \ref{lma4pre}.
\end{Assumption}

    The function $\widehat{V}$ can be solved in the explicit form given in \eqref{vtxz} later. The process  $\widetilde{V}(t,\mu_t; x_0-\kappa t, z_0)$ with the function $\widetilde{V}$ defined in \eqref{tilde2v} is the Snell envelope of the process $\widehat{V}(t, x_0-\kappa t, z_0, \mu_t)$ above. The function $\widetilde{V}$ in  \eqref{tilde2v} can therefore be written as
   \begin{align*}
    \widetilde{V}(t,\eta; x_0-\kappa t, z_0)=\underset{\tau\geq t}{\textrm{esssup}}\ \mathbb{E}\Big[
    \widehat{V}(\tau, x_0-\kappa\tau, z_0, \mu_\tau)
    \Big{|}\mu_t=\eta
    \Big].
    \end{align*}
    The continuation region, interpreted as the continuation of full information observations to update the input value, is denoted by $\mathcal{C}=\{(t,\eta)\in[0,T)\times\mathbb{R}:\widetilde{V}(t,\eta; x_0-\kappa t, z_0)>\widehat{V}(t, x_0-\kappa t,z_0, \eta)\}$ and the free boundary is $\partial\mathcal{C}=\{(t,\eta)\in[0,T)\times\mathbb{R}:\widetilde{V}(t,\eta; x_0-\kappa t, z_0)=\widehat{V}(t, x_0-\kappa t,z_0, \eta)\}$. Let us denote $\widetilde{V}(t,\eta; x_0-\kappa t, z_0)$ by $\widetilde{V}(t,\eta)$ for short when there is no confusion. By some heuristic arguments, we can write the HJB variational inequalities with the terminal condition $\widetilde{V}(T,\eta)=0$, $\eta\in\mathbb{R}$, by
    \begin{align}\label{varine}
    \mbox{min}\left\{ \widetilde{V}(t,\eta) - \widehat{V}(t, x_0-\kappa t, z_0, \eta), \ \
               -\frac{\partial \widetilde{V}(t,\eta)}{\partial t}-\mathcal{L}\widetilde{V}(t,\eta)
       \right\}=0,
    \end{align}
where $\mathcal{L}\widetilde{V}(t,\eta)=
-\lambda(\eta-\bar{\mu})\frac{ \partial \widetilde{V} }{ \partial \eta }(t,\eta)
  +\frac{1}{2} \sigma_\mu^2 \frac{ \partial^2 \widetilde{V} }{ \partial \eta^2 }(t,\eta)$. To simplify notations in the following sections, we shall rewrite \eqref{varine} by
    \begin{align}\label{pde1}
    \left\{\begin{array}{lr}
       F(t,\eta,\widetilde{V}, \frac{\partial\widetilde{V}}{\partial t},\frac{ \partial \widetilde{V} }{ \partial \eta }, \frac{ \partial^2 \widetilde{V} }{ \partial \eta^2 } )=0, \ \ \  \mbox{on} \ [0,T)\times\mathbb{R},\\
       v(T,\eta)=0, \ \ \  \mbox{for} \ \eta\in\mathbb{R},
    \end{array}\right.
    \end{align}
with the operator $F(t,\eta,v,v_t,v_\eta,v_{\eta\eta}):=\mbox{min}\left\{ v-\widehat{V}, \ \
           -\frac{\partial v}{\partial t}-\mathcal{L}v\right\}$.

\begin{remark}
The part $-\frac{\partial \widetilde{V}}{\partial t}-\mathcal{L}\widetilde{V}=0$ in (\ref{varine}) is a linear parabolic PDE and does not depend on the interior control $(\pi,c)$. The comparison part $\widetilde{V} - \widehat{V}$ in (\ref{varine}) depends on the optimal control $(\pi,c)$ because the $\widehat{V}$ is the value function of the interior control problem provided the input $\hat{X}_t=x_0-\kappa t$, $Z_t=z_0$ and $\hat{\mu}_t=\mu_t=\eta$.
\end{remark}

The next theorem is the main result of this paper, whose proof is given in Section \ref{sec-4}.
\begin{theorem}\label{unique}
$\widetilde{V}(t,\eta)$ defined in (\ref{tilde2v}) is the unique bounded and continuous viscosity solution to variational inequalities (\ref{varine}). In addition, the optimal entry time for the composite problem  (\ref{tilde2v}) is given by the $\mathcal{F}_t$-adapted stopping time
\begin{align}\label{optimalentry}
\tau^*:=T\wedge\inf\left\{t\geq 0: \widetilde{V}(t,\mu_t; x_0-\kappa t, z_0)=\widehat{V}(t, x_0-\kappa t,z_0, \mu_t)\right\}.
\end{align}
We also have that the process $\widetilde{V}(t,\mu_t; x_0-\kappa t,z_0)$ is a martingale with respect to the full information filtration $\mathcal{F}_t$ for $0\leq t\leq \tau^*$. \end{theorem}

\section{Interior Utility Maximization under Partial Observations}\label{sec-3}

We first solve the interior stochastic control problem under partial observations of stock prices.

\subsection{Optimal consumption with Kalman-Bucy filtering}  For some fixed time $0\leq k\leq T$, the dynamic interior stochastic control problem under habit formation is defined by
    \begin{align}\label{vxz}
    \begin{aligned}
      &\widehat{V}(k,x,z,\eta;\theta):=\sup_{(\pi,c)\in\mathcal{A}_k(x)}\mathbb{E}\left[\int^{T}_{k}\frac{(c_{s}-Z_{s})^p}{p}ds\Big{|}\mathcal{F}_k^S\right]\\
     =&\sup_{(\pi,c)\in\mathcal{A}_k(x)}\mathbb{E}\left[\int^{T}_{k}\frac{(c_{s}-Z_{s})^p}{p}ds\Big{|}\hat{X}_{k}=x,Z_{k}=z,\hat{\mu}_{k}=\eta;\hat{\Sigma}(k)=\theta\right],
    \end{aligned}
    \end{align}
where $\mathcal{A}_k(x)$ denotes the admissible control space starting from time $k$. Here, as the conditional variance $\hat{\Sigma}(t)$ is a deterministic function of time, we set $\theta$ as a parameter instead of a state variable.

By using the optimality principle and It\^o's formula, we can heuristically obtain the HJB equation as
    \begin{align}\label{hjb}
    \begin{split}
      &V_{t}-\alpha(t)zV_{z}-\lambda(\eta-\bar{\mu})V_{\eta}+\frac{\left(\hat{\Sigma}(t)+\sigma_{S}\sigma_{\mu}\rho\right)^2}{2\sigma_{S}^2}V_{\eta\eta}+\max_{(\pi,c)\in\mathcal{A}}\left[-cV_{x}+c\delta(t)V_{z}+\frac{(c-z)^p}{p}\right]\\
      &+\max_{(\pi,c)\in\mathcal{A}}\left[\pi\eta V_{x}+\frac{1}{2}\sigma_{S}^2\pi^2V_{xx}+V_{x\eta}\left(\hat{\Sigma}(t)+\sigma_{S}\sigma_{\mu}\rho\right)\pi\right]=0,\ \ \ k\leq t\leq T,
    \end{split}
    \end{align}
    with the terminal condition $V(T,x,z,\eta)=0$.

\subsection{The decoupled solution and main results} If $V(t,x,z,\eta)$ is smooth enough, the first order condition gives
    \begin{align*}
    \begin{split}
      & \pi^*(t,x,z,\eta)=\frac{-\eta V_{x}-\left(\hat{\Sigma}(t)+\sigma_{S}\sigma_{\mu}\rho\right)V_{x\eta}}{\sigma_{S}^2 V_{xx}},\\
      & c^*(t,x,z,\eta)=z+\Big(V_{x}-\delta(t)V_{z}\Big)^{\frac{1}{p-1}}.
    \end{split}
    \end{align*}
 Thanks to the homogeneity property of the power utility, we conjecture the value function in the form
    $V(t,x,z,\eta)=\frac{[(x-m(t,\eta)z)]^p}{p}N^{1-p}(t,\eta)$ for some functions $m(t,\eta)$ and $N(t,\eta)$ to be determined. It also follows that the terminal condition that $N(T,\eta)=0$ is required. In particular, we find that the simple ansatz of $m(t,\eta):=m(t)$ satisfies the equation \eqref{mt}. After substitution, the HJB equation reduces to the linear parabolic PDE for $N(t,\eta)$ as
    \begin{align*}
    \begin{aligned}
      N_t&+\frac{p\eta^2}{2(1-p)^2\sigma_S^2}N(t,\eta)
      +\frac{\left(\hat{\Sigma}(t)+\sigma_{S}\sigma_{\mu}\rho\right)^2}{2\sigma_S^2}N_{\eta\eta}
      +\Big(1+\delta(t)m(t)\Big)^{\frac{p}{p-1}}\\
      &+\left[-\lambda(\eta-\bar{\mu})+\frac{\eta\Big(\hat{\Sigma}(t)+\sigma_{S}\sigma_{\mu}\rho\Big)p}{(1-p)\sigma_S^2}\right]N_{\eta}(t,\eta)=0,
    \end{aligned}
    \end{align*}
    with $N(T,\eta)=0$. We can further solve the linear PDE explicitly by
    \begin{align}\label{nteta}
      N(t,\eta)=\int^T_t\Big(1+\delta(s)m(s)\Big)^{\frac{p}{p-1}}\exp\left(A(t,s)\eta^2+B(t,s)\eta+C(t,s)\right)ds,
    \end{align}
    for $k\leq t\leq s\leq T$. $A(t,s)$, $B(t,s)$ and $C(t,s)$ satisfy the following ODEs:
    \begin{align}\label{atts}
     A_t(t,s)+\frac{p}{2(1-p)^2\sigma_S^2}+2\left[-\lambda+\frac{p(\hat{\Sigma}(t)+\sigma_{S}\sigma_{\mu}\rho)}{\sigma_S^2(1-p)}\right]A(t,s)+\frac{2(\hat{\Sigma}(t)+\sigma_{S}\sigma_{\mu}\rho)^2}{\sigma_{S}^2}A^2(t,s)=0,
    \end{align}
    \begin{align}\label{btts}
     &B_t(t,s)+\left[-\lambda+\frac{p(\hat{\Sigma}(t)+\sigma_{S}\sigma_{\mu}\rho)}{\sigma_S^2(1-p)}\right]B(t,s)+2\lambda\bar{\mu}A(t,s)\notag\\
     +&\frac{2(\hat{\Sigma}(t)+\sigma_{S}\sigma_{\mu}\rho)^2}{\sigma_{S}^2}A(t,s)B(t,s)=0,
    \end{align}
    \begin{align}\label{ctts}
      C_t(t,s)+\lambda\bar{\mu}B(t,s)+\frac{\Big(\hat{\Sigma}(t)+\sigma_{S}\sigma_{\mu}\rho\Big)^2}{2\sigma_{S}^2}\Big(B^2(t,s)+2A(t,s)\Big)=0,
    \end{align}
with terminal conditions $A(s,s)=B(s,s)=C(s,s)=0$. The explicit solutions of ODEs (\ref{atts}), (\ref{btts}), (\ref{ctts}) are reported in Appendix \ref{appendix5.1}. For fixed $t\in[k,T]$, we can define the \emph{effective domain} of the pair $(x,z)$ by $\mathbb{D}_t\ :=\ \{(x',z')\in(0,+\infty)\times[0,+\infty);\ x'\geq m(t)z'\}$, where $k\leq t\leq T$. The HJB equation (\ref{hjb}) admits a classical solution on $[k,T]\times\mathbb{D}_t\times\mathbb{R}$ that
     \begin{align}\label{vtxz}
     \begin{aligned}
      &V(t,x,z,\eta)\\
      &=\Big[\int^T_t\Big(1+\delta(s)m(s)\Big)^{\frac{p}{p-1}}\exp\Big(A(t,s)\eta^2+B(t,s)\eta+C(t,s)\Big)ds\Big]^{1-p}\times\frac{[(x-m(t)z)]^p}{p}.
     \end{aligned}
     \end{align}

\begin{remark}\label{effdomain}
The effective domain of $V(t,x,z,\eta)$ requires some constraints on the optimal wealth process $\hat{X}^*_t$ and habit formation process $Z^*_t$ such that $\hat{X}^*_t\geq m(t)Z^*_t$ for $t\in[k,T]$. In particular, we have to enforce the initial wealth-habit budget constraint that $\hat{X}_k\geq m(k)Z_k$ at the initial time $k$.
\end{remark}

\begin{theorem}\label{thm}[The Verification Theorem] If the initial budget constraint $\hat{X}_k\geq m(k)Z_k$ holds at time $k$, the unique solution (\ref{vtxz}) of HJB equation equals the value function defined in (\ref{vxz}), i.e., $V(k,x,z,\eta)=\widehat{V}(k,x,z,\eta)$.
Moreover, the optimal investment policy $\pi^*_t$ and optimal consumption policy $c^*_t$ are given in the feedback form by  $\pi^*_t=\pi^*(t,\hat{X}^*_t,Z^*_t,\hat{\mu}_t)$ and  $\ c^*_t=c^*(t,\hat{X}^*_t,Z^*_t,\hat{\mu}_t)$,  $k\leq t\leq T$. The function $\pi^*(t,x,z,\eta):[k,T]\times\mathbb{D}_t\times\mathbb{R}\rightarrow\mathbb{R}$ is given by
  \begin{align}\label{pistar}
    \pi^*(t,x,z,\eta)
    =\left[\frac{\eta}{(1-p)\sigma^2_S}+\frac{\left(\hat{\Sigma}(t)+\sigma_S\sigma_\mu\rho\right)}{\sigma^2_S}\frac{N_{\eta}(t,\eta)}{N(t,\eta)}\right](x-m(t)z),
  \end{align}
and the function $c^*(t,x,z,\eta):[k,T]\times\mathbb{D}_t\times\mathbb{R}\rightarrow\mathbb{R}^+$ is given by
  \begin{align}\label{cstar}
    c^*(t,x,z,\eta)=z+\frac{(x-m(t)z)}{\Big(1+\delta(t)m(t)\Big)^{\frac{1}{1-p}}N(t,\eta)}.
  \end{align}
\indent The optimal wealth process $\hat{X}^*_t$ , $k\leq t\leq T$, is given by
  \begin{align}\label{xstart}
    \hat{X}^*_t=&(x-m(k)z)\frac{N(t,\hat{\mu}_t)}{N(k,\eta)}\exp\left(\int^t_k\frac{(\hat{\mu}_u)^2}{2(1-p)\sigma^2_S}du+\int^t_k\frac{\hat{\mu}_u}{(1-p)\sigma_S}d\hat{W}_u\right)
    +m(t)Z^*_t.
  \end{align}
\end{theorem}

\section{Exterior Optimal Stopping Problem}\label{sec-4}

\subsection{Stochastic Perron's method}
    We next study the exterior optimal entry problem. Recall that $\hat{X}_{\tau}=x_0-\kappa\tau$, $Z_{\tau}=z_0$ and $\hat{\Sigma}(\tau)=0$ are all taken as parameters. Our aim is to solve an optimal stopping problem in which $\mu_t$ is the only underlying state process.

\begin{remark}\label{convexconcave}
Recall that the interior value function $\widehat{V}$ is of the form in \eqref{vtxz}. Moreover, by Remark  \ref{appendixrk1}, functions $A(t,s)\leq 0$ and $B(t,s)\leq 0$ in \eqref{vtxz} due to $p<0$. That is, if we take $\widehat{V}(\tau, \hat{\mu}_{\tau})$ as a functional of the input $\hat{\mu}_{\tau}$, it is not globally convex or concave in $\hat{\mu}_{\tau}\in\mathbb{R}$ because the function $\exp\left({A(t,s)\eta^2+B(t,s)\eta+C(t,s)}\right)$ is not globally convex or concave in the variable $\eta\in\mathbb{R}$, which depends on values of $A(t,s)$ and $B(t,s)$. Therefore, the composite value function $\widetilde{V}(t,\eta)$ in \eqref{tilde2v} is not globally convex or concave in $\eta\in\mathbb{R}$, which actually depends on all model parameters.
\end{remark}

   We choose to apply the stochastic Perron's method in the present paper to verify that the value function of the composite problem is the unique viscosity solution of some variational inequalities. We first introduce sets of stochastic semi-solutions $\mathcal{V}^+$ and $\mathcal{V}^-$ and prove that $v^-\leq \widetilde{V}\leq v^+$, where $v^-$ and $v^+$ are defined later in (\ref{v-}) and (\ref{v+}). By using the stochastic Perron's method, we can show that $v^+$ is a bounded and upper semi-continuous (u.s.c.) viscosity subsolution and $v^-$ is a bounded and lower semi-continuous (l.s.c.) viscosity supersolution. At last, we prove the comparison principle, that is, if we have any bounded and u.s.c. viscosity subsolution $u$ and bounded and l.s.c. viscosity supersolution $v$ of (\ref{pde1}), we must have the order $u\leq v$. It follows that $v^+\leq v^-$, which leads to the desired conclusion that $v^-= \widetilde{V}= v^+$ and the value function is the unique viscosity solution.

We next present the definitions of stochastic semi-solutions, which are mainly motivated by \cite{bayraktar2014stochastic}.

\begin{definition}\label{def1}
The set of stochastic super-solutions for the PDE (\ref{pde1}), denoted by $\mathcal{V}^+$, is the set of functions $v:[0,T]\times\mathbb{R}\longrightarrow\mathbb{R}$ which have the following properties:\\
    (i) $v$ is u.s.c. and bounded on $[0,T]\times\mathbb{R}$ and $v(t,\eta)\geq \widehat{V}(t,x_0-\kappa t, z_0,\eta)$ for any $(t,\eta)\in[0,T]\times\mathbb{R}$.\\
    (ii) for each $(t,\eta)\in[0,T]\times\mathbb{R}$ and any stopping time $t\leq \tau_1\in \mathcal {T}$, we have $v(\tau_1,\mu_{\tau_1})\geq\mathbb{E}[v(\tau_2,\mu_{\tau_2})|\mathcal{F}_{\tau_1}]$, $\mathbb{P}$-a.s., for any $\tau_2\in\mathcal{T}$ and $\tau_2\geq \tau_1$. That is, the function $v$ along the solution of the SDE (\ref{realmu}) is a supermartingale under the full information filtration $(\mathcal {F}_t)_{t\in[0,T]}$ between $\tau_1$ and $T$.
\end{definition}

\begin{definition}\label{def2}
The set of stochastic sub-solutions for the PDE (\ref{pde1}), denoted by $\mathcal{V}^-$, is the set of functions $v:[0,T]\times \mathbb{R}\longrightarrow\mathbb{R}$ which have the following properties:\\
    (i) $v$ is l.s.c. and bounded on $[0,T]\times \mathbb{R}$ and $v(T,\eta)\leq 0$ for any $\eta\in\mathbb{R}$.\\
    (ii) for each $(t,\eta)\in[0,T]\times \mathbb{R}$ and any stopping time $t\leq \tau_1\in \mathcal {T}$, we have $v(\tau_1,\mu_{\tau_1})\leq\mathbb{E}[v(\tau_2\wedge \zeta,\mu_{\tau_2\wedge \zeta})|\mathcal{F}_{\tau_1}]$, $\mathbb{P}$-a.s., for any $\tau_2\in\mathcal{T}$ and $\tau_2\geq\tau_1$. That is, the function $v$ along the solution to (\ref{realmu}) is a submartingale under the full information filtration $(\mathcal {F}_t)_{t\in[0,T]}$ between $\tau_1$ and $\zeta$, where
\begin{align}\label{hittime}
\zeta:=\inf\{t\in [\tau_1, T]: v(t,\mu_t; x_0-\kappa t, z_0)\geq \widehat{V}(t, x_0-\kappa t,z_0, \mu_t)\}.
\end{align}
\end{definition}

\begin{remark}
We note that the definitions of stochastic super-solutions and stochastic sub-solutions for the optimal stopping problem are not symmetric, which are consistent with the similar definitions in \cite{bayraktar2014stochastic}. The main reason for these differences comes from the natural supermartingale property of the Snell envelop process and its martingale property between the initial time and the first hitting time $\zeta$ in \eqref{hittime}. That is, we naturally need $v(t,\eta)\geq \widehat{V}(t,x_0-\kappa t, z_0,\eta)$ for all $(t,\eta)\in[0,T]\times\mathbb{R}$ including the terminal time $T$ in item $(i)$ of Definition \ref{def1} of stochastic super-solution, but we only require $v(T,\eta)\leq \widehat{V}(T,x_0-\kappa t, z_0,\eta)=0$ at the terminal time $T$ in item $(i)$ of Definition \ref{def2} for stochastic sub-solution. These comparison results and the supermartingale and submartingale properties will play important roles to establish the desired sandwich result $v^-\leq \widetilde{V}\leq v^+$ in Lemma \ref{compalem}.
\end{remark}

\begin{lemma}\label{bddcon}
$\widehat{V}(t,x_0-\kappa t, z_0, \eta; 0)$ is bounded and continuous for $(t,\eta)\in [0,T]\times \mathbb{R}$.
\end{lemma}
\begin{proof} For fixed $x_0$ and $z_0$, it is clear that $\widehat{V}(t,x_0-\kappa t,z_0,\eta)$ in \eqref{vtxz} is continuous and $\widehat{V}(t,x_0-\kappa t,z_0,\eta)\leq 0$. Therefore we only need to show that $\widehat{V}$ is lower bounded. By Appendix \ref{appendix5.1}, we know that $A(u)\leq 0$, $B(u)\leq 0$ and $C(u)\leq K$ for some $K\geq 0$ thanks to $p<0$. We hence obtain that $\Big(A(u)\eta^2+B(u)\eta+C(u)\Big)\leq K_1$ for some $K_1>0$ and it follows that $\widehat{V}(t,x_0-\kappa t,z_0,\eta)$ is lower bounded by some constant for $(t,\eta)\in [0,T]\times \mathbb{R}$ again by $p<0$.
\end{proof}

As it is trivial to see that $0\in\mathcal{V}^-$ and $0\in\mathcal{V}^+$, we have the following result.

\begin{lemma}\label{nonempty}
$\mathcal{V}^+$ and $\mathcal{V}^-$ are nonempty.
\end{lemma}

\begin{definition}\label{v-+}
   We define
   \begin{align}
   v^-:=\sup_{p\in\mathcal{V}^-}p\label{v-},\\ v^+:=\inf_{q\in\mathcal{V}^+}q\label{v+}.
   \end{align}
\end{definition}

Similar to Lemma 2.2. of \cite{bayraktar2011probabilistic}, the next result holds.

\begin{lemma}\label{in}
We have $v^-\in\mathcal{V}^-$ and $v^+\in\mathcal{V}^+$.
\end{lemma}

We have the first important sandwich result.

\begin{lemma}\label{compalem}
We have $v^-\leq \widetilde{V}\leq v^+$.
\end{lemma}

\begin{proof}  For each $v\in\mathcal{V}^+$, let us consider $\tau_1=t\geq 0$ in Definition \ref{def1}. For any $\tau\geq t$, we have $v(t,\eta)\geq\mathbb{E}[v(\tau,\mu_{\tau})|\mathcal{F}_t]\geq\mathbb{E}[\widehat{V}(\tau,x_0-\kappa\tau,z_0,\mu_\tau)|\mathcal{F}_t]$ thanks to the supermartingale property in Definition \ref{def1}. It follows that $v(t,\eta)\geq{\textrm{esssup}_{t\leq \tau}}\ \mathbb{E}[\widehat{V}(\tau,x_0-\kappa\tau,z_0, \mu_\tau)|\mathcal{F}_t]$. This implies that $v(t,\eta)\geq\widetilde{V}(t,\eta)$ in view of the definition of $\widetilde{V}(t,\eta)$ and hence $\widetilde{V}\leq v^+$ by the definition in (\ref{v+}). On the other hand, for each $v\in\mathcal{V}^-$, by taking $\tau_1=t\geq 0$ in Definition \ref{def2}, we have $v(t,\eta)\leq\mathbb{E}[v(\tau\wedge\zeta,\mu_{\tau\wedge\zeta})|\mathcal{F}_t]$ for any $\tau\geq t$ because of the submartingale property in Definition \ref{def2}. In particular, using the definition of $\zeta$, we further have $v(t,\eta)\leq\mathbb{E}[v(\tau\wedge\zeta,\mu_{\tau\wedge\zeta})|\mathcal{F}_t]\leq \mathbb{E}[\widehat{V}(\tau\wedge\zeta,x_0-f(\tau\wedge\zeta),z_0,\mu_{\tau\wedge\zeta})|\mathcal{F}_t]\leq{\textrm{esssup}_{\tau\geq t}}\ \mathbb{E}[\widehat{V}(\tau,x_0-\kappa\tau,z_0, \mu_\tau)|\mathcal{F}_t]=\widetilde{V}(t,\eta)$. It then follows that $\widetilde{V}\geq v^-$ because of (\ref{v-}). In conclusion, we have the inequality $v^-\leq \widetilde{V}\leq v^+$.
\end{proof}

\begin{theorem}\label{SPM} $v^-$ in Definition \ref{v-+} is a bounded and l.s.c. viscosity super-solution of
    \begin{align}\label{vissupsol}
    \left\{\begin{array}{lr}
       F(t,\eta,v,v_t,v_\eta,v_{\eta\eta})\geq 0, \ \ \  \mbox{on} \ [0,T)\times\mathbb{R},\\
       v(T,\eta)\geq 0, \ \ \  \mbox{for any} \ \eta\in\mathbb{R},
      \end{array}\right.
    \end{align}
and $v^+$ in Definition \ref{v-+} is a bounded and u.s.c. viscosity sub-solution of
    \begin{align}\label{vissubsol}
    \left\{\begin{array}{lr}
       F(t,\eta,v,v_t,v_\eta,v_{\eta\eta})\leq 0, \ \ \  \mbox{on} \ [0,T)\times\mathbb{R},\\
       v(T,\eta)\leq 0, \ \ \  \mbox{for any} \ \eta\in\mathbb{R}.
      \end{array}\right.
    \end{align}
\end{theorem}

\begin{proof} We follow and modify some arguments in \cite{bayraktar2011probabilistic,bayraktar2014stochastic} to fit our setting.

\textit{(i) The sub-solution property of $v^+$.} First, definition in \eqref{v+} and Lemma \ref{in} imply that $v^+$ is bounded and upper semi-continuous. Suppose $v^+$ is not a viscosity sub-solution, there exists some interior point $(\bar{t},\bar{\eta})\in(0,T)\times\mathbb{R}$ and a $C^{1,2}$-test function $\varphi:[0,T]\times\mathbb{R}\rightarrow\mathbb{R}$ such that $v^+ - \varphi$ attains a strict local maximum that is equal to zero and $F(\bar{t},\bar{\eta},v,v_{\bar{t}},v_{\bar{\eta}},v_{\bar{\eta}\bar{\eta}})>0$. It follows that
    \begin{align*}
    \left\{
       \begin{array}{lr}
       v^+(\bar{t},\bar{\eta}) - \widehat{V}(\bar{t},x_0-f(\bar{t}),z_0,\bar{\eta})> 0,\\
       -\frac{\partial \varphi}{\partial t}(\bar{t},\bar{\eta})-\mathcal{L}\varphi(\bar{t},\bar{\eta}) >0.
       \end{array}
    \right.
    \end{align*}
  There exists a ball $B(\bar{t},\bar{\eta},\varepsilon)$ small enough that
  \begin{align*}
  \left\{
  \begin{array}{lr}
  -\frac{\partial \varphi}{\partial t}-\mathcal{L}\varphi>0 \ \ \mbox{on} \ \ \overline{B(\bar{t},\bar{\eta},\varepsilon)},\\
  \varphi>v^+ \ \ \mbox{on} \ \ \overline{B(\bar{t},\bar{\eta},\varepsilon)}\backslash(\bar{t},\bar{\eta}).
  \end{array}
  \right.
  \end{align*}

  In addition, as $\varphi(\bar{t},\bar{\eta}) = v^+(\bar{t},\bar{\eta})> \widehat{V}(\bar{t},x_0-f(\bar{t}),z_0,\bar{\eta})$, $\varphi$ is continuous and $\widehat{V}$ is continuous, we can derive that for some $\varepsilon$ small enough, we have  $\varphi-\varepsilon\geq\widehat{V}$ on $\overline{B(\bar{t},\bar{\eta},\varepsilon)}$. Because $v^+-\varphi$ is upper semi-continuous and $\overline{B(\bar{t},\bar{\eta},\varepsilon)}\backslash B(\bar{t},\bar{\eta},\frac{\varepsilon}{2})$ is compact, it then follows that there exists a $\delta>0$ such that $\varphi-\delta\geq v^+$ on $\overline{B(\bar{t},\bar{\eta},\varepsilon)}\backslash B(\bar{t},\bar{\eta},\frac{\varepsilon}{2})$.

  If we choose $0<\xi<\delta\wedge\varepsilon$, the function $\varphi^\xi=\varphi-\xi$ satisfies that
  \begin{align*}
  \left\{
  \begin{array}{lr}
  -\frac{\partial \varphi^\xi}{\partial t}-\mathcal{L}\varphi^\xi>0 \ \  \mbox{on} \ \ \overline{B(\bar{t},\bar{\eta},\varepsilon)},\\
  \varphi^\xi>v^+ \ \ \mbox{on} \ \  \overline{B(\bar{t},\bar{\eta},\varepsilon)}\backslash B(\bar{t},\bar{\eta},\frac{\varepsilon}{2}),\\
  \varphi^\xi\geq\widehat{V} \ \  \mbox{on} \ \ \overline{B(\bar{t},\bar{\eta},\varepsilon)},
  \end{array}
  \right.
  \end{align*}
and $\varphi^\xi(\bar{t},\bar{\eta})=v^+(\bar{t},\bar{\eta})-\xi$.

Let us define an auxiliary function by
  \begin{align*}
  v^\xi:=\left\{
  \begin{array}{lr}
  v^+\wedge\varphi^\xi \ \  \mbox{on} \ \ \overline{B(\bar{t},\bar{\eta},\varepsilon)},\\
  v^+ \ \ \mbox{outside} \ \ \overline{B(\bar{t},\bar{\eta},\varepsilon)}.
  \end{array}
  \right.
  \end{align*}
It is easy to check that $v^\xi$ is upper semi-continuous and $v^\xi(\bar{t},\bar{\eta})=\varphi^\xi(\bar{t},\bar{\eta})<v^+(\bar{t},\bar{\eta})$. We claim that $v^\xi$ satisfies the terminal condition. To this end, we pick some $\varepsilon>0$ that satisfies $T>\bar{t}+\varepsilon$ and recall that $v^+$ satisfies the terminal condition. We then continue to show that $v^\xi\in \mathcal{V}^+$ to obtain a contradiction.

  Let us fix $(t,\eta)$ and recall that $((\mu_s)_{t\leq s\leq T},(W_s,B_s)_{t\leq s\leq T},\Omega,\mathcal{F},\mathbb{P},(\mathcal{F}_s)_{t\leq s\leq T})\in\chi$, where $\chi$ is the nonempty set of all weak solutions. We need to show that the process $(v^\xi(s,\mu_s))_{t\leq s\leq T}$  is a supermartingale on $(\Omega,\mathbb{P})$ with respect to $(\mathcal{F}_s)_{t\leq s\leq T}$. We first assume that $(v^+(s,\mu_s))_{t\leq s\leq T}$ has right continuous paths. In this case, $v^\xi$ is a supermartingale locally in the region $[t,T]\times\mathbb{R}\backslash B(\bar{t},\bar{\eta},\frac{\varepsilon}{2})$ because it equals the right continuous supermartingale $(v^+(s,\mu_s))_{t\leq s\leq T}$. As the process $(v^\xi(s,\mu_s))_{t\leq s\leq T}$ is the minimum between two local supermartingales in the region $B(\bar{t},\bar{\eta},\varepsilon)$, it is a local supermartingale. As two regions $[t,T]\times\mathbb{R}\backslash B(\bar{t},\bar{\eta},\frac{\varepsilon}{2})$ and $B(\bar{t},\bar{\eta},\varepsilon)$ overlap over an open region, $(v^\xi(s,\mu_s))_{t\leq s\leq T}$ is actually a supermartingale.

 If the process $(v^+(s,\mu_s))_{t\leq s\leq T}$ is not right continuous, we can consider its right continuous limit over rational times to transform it to the special case discussed above. In particular, for a given rational number $r$ and fixed $0 \leq t\leq r \leq s \leq T$ and $\eta\in\mathbb{R}$, it remains to show the process $(Y_u)_{t\leq u\leq T}:=(v^\xi(u,\mu_u))_{t\leq u\leq T}$ between $r$ and $s$ is a supermartingale, which is equivalent to show that $Y_r \geq \mathbb{E}[Y_s|\mathcal{F}_r]$.

  Let us denote $G_u:=v^+(u,\mu_u), \ r \leq u \leq s$ and freeze the process $G$ after time s, i.e. $G_u:=v^+(s,\mu_s), \ s \leq u \leq T$. As $(G_u)_{r\leq u\leq T}$ may not be right continuous, by Proposition 1.3.14 in \cite{motion1988brownian}, we can consider its right continuous modification by
  \begin{align*}
  G^+_u(\omega):=\lim_{u'\rightarrow u, \ u'>u, \ u'\in\mathbb{Q}}G_{u'}(\omega), \ \ \ r\leq u\leq T.
  \end{align*}
Note that $G^+$ is a right continuous supermartingale with respect to $\mathcal{F}$ that satisfies the usual conditions. Because $v^+$ is upper semi-continuous and the process remains the same after $s$, we conclude that $G_r\geq G_r^+, \ G_s=G_s^+$. Recall that $v^+<\varphi-\delta$ in the open region $B(\bar{t},\bar{\eta},\varepsilon)\backslash\overline{B(\bar{t},\bar{\eta},\frac{\varepsilon}{2})}$, if we take right limits inside this region and use continuous function $\varphi$, we have
  \begin{align*}
  G^+_u<\varphi^\xi(u,\mu_u), \ \mbox{if}
  \ (u,\mu_u)\in
  B(\bar{t},\bar{\eta},\varepsilon)\backslash\overline{B(\bar{t},\bar{\eta},\frac{\varepsilon}{2})}.
  \end{align*}
Thus, if we consider the process
  \begin{align*}
  Y_u^+:=\left\{
  \begin{array}{lr}
  G_u^+, \ (u,\mu_u)\not\in\overline{B(\bar{t},\bar{\eta},\frac{\varepsilon}{2})},\\
  G_u^+\wedge\varphi^\xi(u,\mu_u), \ (u,\mu_u)\in B(\bar{t},\bar{\eta},\varepsilon),
  \end{array}
  \right.
  \end{align*}
we also have $Y_r\geq Y_r^+, \ Y_s=Y_s^+$.

Because $G^+$ has right continuous paths, we can conclude that $Y$ is a supermartingale such that
  \begin{align*}
  Y_r\geq Y_r^+\geq\mathbb{E}[Y_s^+|\mathcal{F}_r]
  =\mathbb{E}[Y_s|\mathcal{F}_r].
  \end{align*}

\textit{(ii) The terminal condition of $v^+$.}

For some $\eta_0\in\mathbb{R}$, we assume that $v^+(T,\eta_0)>0$ and will show a contradiction. As $\widehat{V}$ is continuous on $\mathbb{R}$, we can choose an $\varepsilon>0$ such that $0 \leq v^+(T,\eta_0)-\varepsilon$ and $|\eta-\eta_0|\leq\varepsilon$. On the compact set $(\overline{B(T,\eta_0,\varepsilon)}\backslash B(T,\eta_0,\frac{\varepsilon}{2}))\cap([0,T]\times\mathbb{R})$, $v^+$ is bounded above by the definition of $\mathcal{V}^+$ and that $v^+\in\mathcal{V}^+$. Moreover, as $v^+$ is upper semi-continuous on this compact set, we can find $\delta>0$ small enough such that
  \begin{align}\label{initialcon1}
    v^+(T,\eta_0) + \frac{\varepsilon^2}{4\delta} \geq \varepsilon+\sup_{(t,\eta)\in(\overline{B(T,\eta_0,\varepsilon)}\backslash B(T,\eta_0,\frac{\varepsilon}{2}))\cap([0,T]\times\mathbb{R})}v^+(t,\eta).
  \end{align}

Next, for $k>0$, we define the function $\varphi^{\delta,\varepsilon,k}(t,\eta) := v^+(T,\eta_0) + \frac{|\eta-\eta_0|^2}{\delta} + k(T-t)$. For $k$ large enough, we derive that $-\varphi^{\delta,\varepsilon,k}_t - \mathcal{L}\varphi^{\delta,\varepsilon,k} > 0 \ \mbox{on} \  \overline{B(T,\eta_0,\varepsilon)}$. Moreover, we have the following result in view of (\ref{initialcon1})
  \begin{align*}
    \varphi^{\delta,\varepsilon,k} \geq \varepsilon + v^+  \ \mbox{on} \ (\overline{B(T,\eta_0,\varepsilon)}\backslash B(T,\eta_0,\frac{\varepsilon}{2}))\cap([0,T]\times\mathbb{R}),
  \end{align*}
and $\varphi^{\delta,\varepsilon,k}(T,\eta) \geq v^+(T,\eta_0) \geq 0+ \varepsilon$ for $|\eta-\eta_0|\leq\varepsilon$.

Now, we can find $\xi<\varepsilon$ and define the function as follows,
  \begin{align*}
  v^{\delta,\varepsilon,k,\xi}:=\left\{\begin{array}{lr}
  v^+\wedge(\varphi^{\delta,\varepsilon,k}-\xi) \ \  \mbox{on} \ \ \overline{B(T,\eta_0,\varepsilon)},\\
  v^+ \ \ \mbox{outside} \ \ \overline{B(T,\eta_0,\varepsilon)}.
  \end{array}\right.
  \end{align*}
By following similar argument in $(i)$, one can obtain that $v^{\delta,\varepsilon,k,\xi}\in\mathcal{V}^+$, but $v^{\delta,\varepsilon,k,\xi}(T,\eta_0) = v^+(T,\eta_0)-\xi$, which leads to a contradiction.

\textit{(iii) The super-solution property of $v^-$.}

Let us only provide a sketch of the proof as it is essentially similar to Step $(i)$. Suppose that $v^-$ is not a viscosity super-solution, then there exist some interior point $(\bar{t},\bar{\eta})\in(0,T)\times\mathbb{R}$ and a $C^{1,2}$-test function $\psi:[0,T]\times\mathbb{R}\rightarrow\mathbb{R}$ such that $v^- - \psi$ attains a strict local minimum that is equal to zero. As $F(\bar{t},\bar{\eta},v,v_{\bar{t}},v_{\bar{\eta}},v_{\bar{\eta}\bar{\eta}})<0$, there are two separate cases to check.

    case(i) \ $v^-(\bar{t},\bar{\eta})-\widehat{V}(\bar{t},x_0-f(\bar{t}),z_0,\bar{\eta})<0$. This already leads to a contradiction with $v^-(\bar{t},\bar{\eta})\geq\widehat{V}(\bar{t},x_0-f(\bar{t}),z_0,\bar{\eta})$ by the definition of $v^-$.

    case(ii) \ $-\frac{\partial \psi}{\partial t}(\bar{t},\bar{\eta})-\mathcal{L}\psi(\bar{t},\bar{\eta}) <0$. We can find a ball $B(\bar{t},\bar{\eta},\varepsilon)$ small enough such that $-\frac{\partial \psi}{\partial t}-\mathcal{L}\psi<0$ on $\overline{B(\bar{t},\bar{\eta},\varepsilon)}$. Moreover, as $v^--\psi$ is lower semi-continuous and $\overline{B(\bar{t},\bar{\eta},\varepsilon)}\backslash B(\bar{t},\bar{\eta},\frac{\varepsilon}{2})$ is compact, there exists a $\delta>0$ such that $\psi+\delta\leq v^-$ on $\overline{B(\bar{t},\bar{\eta},\varepsilon)}\backslash B(\bar{t},\bar{\eta},\frac{\varepsilon}{2})$. We can then choose $\xi\in(0,\frac{\delta}{2})$ small such that $\psi^\xi=\psi+\xi$ satisfies three properties: (i) $-\frac{\partial \psi^\xi}{\partial t}-\mathcal{L}\psi^\xi<0$ on $\overline{B(\bar{t},\bar{\eta},\varepsilon)}$; (ii) we have $v^-\geq\psi+\delta > \psi+\xi = \psi^\xi$ on $\overline{B(\bar{t},\bar{\eta},\varepsilon)}\backslash B(\bar{t},\bar{\eta},\frac{\varepsilon}{2})$; (iii)  $\psi^\xi(\bar{t},\bar{\eta}) = \psi(\bar{t},\bar{\eta}) + \xi = v^-(\bar{t},\bar{\eta}) +\xi > v^-(\bar{t},\bar{\eta})$. Thus, we can define an auxiliary function by
  \begin{align*}
  v^\xi:=\left\{
  \begin{array}{lr}
  v^-\vee\psi^\xi \ \  \mbox{on} \ \ \overline{B(\bar{t},\bar{\eta},\varepsilon)},\\
  v^- \ \ \mbox{outside} \ \ \overline{B(\bar{t},\bar{\eta},\varepsilon)}.
  \end{array}
  \right.
  \end{align*}

By repeating similar argument in Step $(i)$, we have that $v^\xi\in \mathcal{V}^-$ by showing that $(v^\xi(s,\mu_s))_{t\leq s\leq T}$ is a submartingale. If $v^-$ has right continuous paths, the proof is trivial. In general, by Proposition 1.3.14 in \cite{motion1988brownian}, we can consider the right continuous submartingale $G^+_u(\omega):=\lim_{u'\rightarrow u, \ u'>u, \ u'\in\mathbb{Q}}G_{u'}(\omega), \ \omega\in\Omega^\ast, \ r\leq u\leq T$, where $G_u:=v^-(u, \mu_u), \ r \leq u \leq s$ and we stop it at time $t$. Similar to Step $(i)$, we note that $G^+$ is the right continuous submartingale and therefore $G_r\leq G_r^+, \ G_s=G_s^+$. As $ G^+_u>\psi^\xi(u, \mu_u), \ \mbox{if}
  \ (u,\mu_u)\in
  B(\bar{t},\bar{\eta},\varepsilon)\backslash \overline{B(\bar{t},\bar{\eta},\frac{\varepsilon}{2})}$, we can define the process
  \begin{align*}
  Y_u^+:=\left\{
  \begin{array}{lr}
  G_u^+, \ (u, \mu_u)\not\in\overline{B(\bar{t}, \bar{\eta},\frac{\varepsilon}{2})},\\
  G_u^+\vee\psi^\xi(u, \mu_u), \ (u, \mu_u)\in \overline{B(\bar{t},\bar{\eta},\frac{\varepsilon}{2})}.
  \end{array}
  \right.
  \end{align*}
We can conclude that $Y_r\leq Y_r^+, \ Y_s=Y_s^+$ and $Y$ is a submartingale that $Y_r\leq Y_r^+\leq\mathbb{E}[Y_s^+|\mathcal{F}_r]=\mathbb{E}[Y_s|\mathcal{F}_r]$, which completes the proof.

\textit{(iv) The terminal condition of $v^-$.}

For some $\eta_0\in\mathbb{R}$, suppose that $v^-(T,\eta_0)<0$ and we will show a contradiction. As $\widehat{V}$ is continuous on $\mathbb{R}$, we can choose an $\varepsilon>0$ such that $0 \geq v^-(T,\eta_0)+\varepsilon$ and $|\eta-\eta_0|\leq\varepsilon$. Similar to Step $(ii)$, we can find $\delta>0$ small enough such that
  \begin{align}\label{initialcon2}
    v^-(T,\eta_0) - \frac{\varepsilon^2}{4\delta} \leq \inf_{(t,\eta)\in(\overline{B(T,\eta_0,\varepsilon)}\backslash B(T,\eta_0,\frac{\varepsilon}{2}))\cap([0,T]\times\mathbb{R})}v^-(t,\eta)-\varepsilon.
  \end{align}
Then, for $k>0$, we consider $ \psi^{\delta,\varepsilon,k}(t,\eta) := v^-(T,\eta_0) - \frac{|\eta-\eta_0|^2}{\delta} - k(T-t)$. For $k$ large enough, we have that $-\psi^{\delta,\varepsilon,k}_t - \mathcal{L}\psi^{\delta,\varepsilon,k} < 0 \ \mbox{on} \  \overline{B(T,\eta_0,\varepsilon)}$. Furthermore, in view of (\ref{initialcon2}), we have
  \begin{align*}
    \psi^{\delta,\varepsilon,k} \leq v^- - \varepsilon  \ \mbox{on} \ (\overline{B(T,\eta_0,\varepsilon)}\backslash B(T,\eta_0,\frac{\varepsilon}{2}))\cap([0,T]\times\mathbb{R}),
  \end{align*}
and $\psi^{\delta,\varepsilon,k}(T,\eta) \leq v^-(T,\eta_0) \leq  - \varepsilon$ for $|\eta-\eta_0|\leq\varepsilon$.

Next, we can find $\xi<\varepsilon$ and define the function by
  \begin{align*}
  v^{\delta,\varepsilon,k,\xi}:=\left\{\begin{array}{lr}
  v^-\vee(\psi^{\delta,\varepsilon,k}+\xi) \ \  \mbox{on} \ \ \overline{B(T,\eta_0,\varepsilon)},\\
  v^- \ \ \mbox{outside} \ \ \overline{B(T,\eta_0,\varepsilon)}.
  \end{array}\right.
  \end{align*}
Similar to Step $(iii)$, we obtain that $v^{\delta,\varepsilon,k,\xi}\in\mathcal{V}^-$, but $v^{\delta,\varepsilon,k,\xi}(T,\eta_0) = v^-(T,\eta_0)+\xi$, which gives a contradiction.
\end{proof}

Let us then reverse the time and consider $s:=T-t$. However, for the simplicity of presentation, let us continue to use $t$ in the place of $s$ if there is no confusion. The variational inequalities can be written by
\begin{align}\label{varine2}
    \mbox{min}\left\{ \widetilde{V}(t,\eta;x_0-f(T-t), z_0) - \widehat{V}(t, x_0-f(T-t), z_0, \eta), \ \
               \frac{\partial \widetilde{V}(t,\eta)}{\partial t}-\mathcal{L}\widetilde{V}(t,\eta)
       \right\}=0,
\end{align}
where $\mathcal{L}\widetilde{V}(t,\eta)=
-\lambda(\eta-\bar{\mu})\frac{ \partial \widetilde{V} }{ \partial \eta }(t,\eta)
  +\frac{1}{2} \sigma_\mu^2 \frac{ \partial^2 \widetilde{V} }{ \partial \eta^2 }(t,\eta)$ with the condition $\widetilde{V}(0,\eta)=0$.

  Let us denote it equivalently as
    \begin{align}\label{pde2}
    \left\{\begin{array}{lr}
       F(t,\eta,v,v_t,v_\eta,v_{\eta\eta})=0,  \mbox{on} \ (0,T]\times\mathbb{R},\\
       v(0,\eta)=\widehat{V}(0,x_0-f(0),z_0,\eta),  \mbox{for any} \ \eta\in\mathbb{R},
      \end{array}\right.
    \end{align}
where $F(t,\eta,v,v_t,v_\eta,v_{\eta\eta}):=\mbox{min}\left\{ v-\widehat{V}, \ \
           \frac{\partial v}{\partial t}-\mathcal{L}v\right\}$. We also have the continuation region as $\mathcal{C}=\{(t,\eta)\in(0,T]\times\mathbb{R}:\widetilde{V}(t,\eta;x_0-f(T-t), z_0)>\widehat{V}(t, x_0-f(T-t), z_0, \eta)\}$.

\begin{proposition}\label{compp}[Comparison Principle]
Let $u,v$ be u.s.c viscosity subsolution and l.s.c. viscosity supersolution of (\ref{pde2}), respectively. If $u(0,\eta)\leq v(0,\eta)$ on $\mathbb{R}$, then we have $u\leq v$ on $(0,T]\times\mathbb{R}$.
\end{proposition}

\begin{proof}  We will follow and modify some arguments in \cite{bayraktar2015stochastic,pham2009continuous} to fit our setting. Suppose that $u(0,\eta)\leq v(0,\eta)$ on $\mathbb{R}$, and we will prove that $u\leq v$ on $[0,T]\times\mathbb{R}$. We first construct the strict supersolution to the system (\ref{pde2}) with suitable perturbations of $v$. Let us recall that $A\leq 0$, $B\leq 0$ and $C$ is bounded above by some constant in Remark \ref{appendixrk1}. Moreover, we know that $\widehat{V}(t,x_0-\kappa t,z_0,\eta)\leq 0$. Let us fix a constant $C_2>0$ small enough such that $\lambda>C_2\sigma^2_{\mu}$ and set $\psi(t,\eta) = C_0e^{t}+e^{C_2\eta^2}$ with some $C_0>1$. We have that
  \begin{align*}
  \begin{aligned}
    \frac{\partial \psi}{\partial t} - \mathcal{L}\psi = &C_0e^{t} +C_2\Big[ 2(\lambda-C_2\sigma_\mu^2)\eta^2-2\lambda\bar{\mu}\eta-\sigma_\mu^2 \Big]e^{C_2\eta^2}\\
\geq& C_0e^{t} + C_2\frac{-2(\lambda-C_2\sigma_\mu^2)\sigma_\mu^2-\lambda^2\bar{\mu}^2}{2(\lambda-C_2\sigma_\mu^2)}\\
>& C_0 + C_2\frac{-2(\lambda-C_2\sigma_\mu^2)\sigma_\mu^2-\lambda^2\bar{\mu}^2}{2(\lambda-C_2\sigma_\mu^2)}.
  \end{aligned}
  \end{align*}
We can then choose $C_0>1$ large enough such that $C_0 + C_2\frac{-2(\lambda-C_2\sigma_\mu^2)\sigma_\mu^2-\lambda^2\bar{\mu}^2}{2(\lambda-C_2\sigma_\mu^2)}>1$, which guarantees that
\begin{align}\label{A}
\frac{\partial \psi}{\partial t} - \mathcal{L}\psi >1.
\end{align}

Let us define $v^\Lambda := (1-\Lambda)v + \Lambda\psi$ on $[0,T]\times\mathbb{R}$ for any $\Lambda\in(0,1)$. It follows that
     \begin{align}\label{B}
     \begin{split}
        v^\Lambda - \widehat{V} &= (1-\Lambda)v + \Lambda\psi - \widehat{V}= (1-\Lambda)v + \Lambda (C_0e^{t}+e^{C_2\eta^2})  -\widehat{V}\\
      & \geq (1-\Lambda)v + \Lambda (C_0e^{t}+e^{C_2\eta^2}) + \Lambda\widehat{V} -\widehat{V}\\
      &> (1-\Lambda)(v - \widehat{V}) + \Lambda C_0 > \Lambda,
     \end{split}
     \end{align}
    where we used $v - \widehat{V} \geq 0$ in the last inequality. From (\ref{A}) and (\ref{B}), we can deduce that for $\Lambda\in(0,1)$, $v^\Lambda$ is a supersolution to
     \begin{align}\label{C}
      \mbox{min}\left\{ v^\Lambda-\widehat{V}, \ \
           \frac{\partial v^\Lambda}{\partial t}-\mathcal{L}v^\Lambda\right\}
           \geq \Lambda.
     \end{align}

    In order to prove the comparison principle, it suffices to show the claim that $\sup(u-v^\Lambda)\leq 0$ for all $\Lambda\in(0,1)$, as the required result is obtained by letting $\Lambda$ go to $0$. To this end, we will prove the claim by showing a contradiction and suppose that there exists some $\Lambda\in(0,1)$ such that $M:=\sup(u-v^\Lambda)>0$.

    It is clear that $u$, $v$ and $\widehat{V}$ have the same growth conditions: in view of the explicit forms of $A,B,C$ and $\widehat{V}$, it follows that $\widehat{V}$ has growth condition in $t$ as $e^{e^{K_1t}}$ for some $K_1<0$ and has growth condition in $\eta$ as $e^{K_2\eta^2}$ for some $K_2<0$; on the other hand, $\psi$ has growth condition in $t$ as $e^{t}$ and has growth condition in $\eta$ as $e^{C_2\eta^2}$. Thus, we have that $u(t,\eta)-v^\Lambda(t,\eta) = ( u - (1-\Lambda)v - \Lambda\psi)(t,\eta)$ goes to $-\infty$ as $t\rightarrow T,\eta\rightarrow\infty$. Consequently, the u.s.c. function $(u-v^\Lambda)$ attains its maximum $M$.

    Let us consider the u.s.c. function $\Phi_\varepsilon(t,t',\eta,\eta')=u(t,\eta)-v^\Lambda(t',\eta')-\phi_\varepsilon(t,t',\eta,\eta')$, where $\phi_\varepsilon(t,t',\eta,\eta') = \frac{1}{2\varepsilon}((t-t')^2+(\eta-\eta')^2)$, $\varepsilon >0$ and $(t_\varepsilon,t'_\varepsilon,\eta_\varepsilon,\eta'_\varepsilon)$ attains the maximum of $\Phi_\varepsilon$. We have
     \begin{align}\label{M}
     M_\varepsilon = \max\Phi_\varepsilon = \Phi_\varepsilon(t_\varepsilon,t'_\varepsilon,\eta_\varepsilon,\eta'_\varepsilon)\rightarrow
     M \ \text{and} \ \phi_\varepsilon(t_\varepsilon,t'_\varepsilon,\eta_\varepsilon,\eta'_\varepsilon)\rightarrow 0 \ \text{when} \ \varepsilon\rightarrow 0.
     \end{align}

    Let us recall the equivalent definition of viscosity solutions in terms of superjets and subjets. In particular, we define $\bar{\mathcal{P}}^{2,+}u(\bar{t},\bar{\eta})$ as the set of elements $(\bar{q},\bar{k},\bar{M})\in\mathbb{R}\times\mathbb{R}\times\mathbb{R}$ satisfying $u(t,\eta)\leq u(\bar{t},\bar{\eta}) + \bar{q}(t-\bar{t})+ \bar{k}(\eta-\bar{\eta})+ \frac{1}{2}\bar{M}(\eta-\bar{\eta})^2 +  o((t-\bar{t})+(\eta-\bar{\eta})^2)$. We define $\bar{\mathcal{P}}^{2,-}v^\Lambda(\bar{t},\bar{\eta})$ similarly. Thanks to Crandall-Ishii's lemma, we can find $A_\varepsilon, B_\varepsilon\in\mathbb{R}$ such that
     \begin{align*}
     \begin{split}
     (\frac{t_\varepsilon-t_\varepsilon'}{\varepsilon},\frac{\eta_\varepsilon-\eta_\varepsilon'}{\varepsilon},A_\varepsilon) &\in \bar{\mathcal{P}}^{2,+}u(t_\varepsilon,\eta_\varepsilon), \\
     (\frac{t_\varepsilon-t_\varepsilon'}{\varepsilon},\frac{\eta_\varepsilon-\eta_\varepsilon'}{\varepsilon},B_\varepsilon) &\in \bar{\mathcal{P}}^{2,-}v^\Lambda(t_\varepsilon',\eta_\varepsilon'), \\
     \sigma^2(\eta_\varepsilon)A_{\varepsilon} - \sigma^2(\eta_\varepsilon')B_{\varepsilon} &\leq \frac{3}{\varepsilon}(\sigma(\eta_\varepsilon)-\sigma(\eta_\varepsilon'))^2.
     \end{split}
     \end{align*}

    By combining the viscosity subsolution property (\ref{vissubsol}) of $u$ and the viscosity strict supersolution property (\ref{C}) of $v^\Lambda$, we have that
     \begin{align}
     \begin{split}
      \mbox{min}\Big\{ u(t_\varepsilon,\eta_\varepsilon)-\widehat{V}(t_\varepsilon,x_0-f(t_\varepsilon),z_0,\eta_\varepsilon),\frac{t_\varepsilon-t_\varepsilon'}{\varepsilon}
           -\frac{\eta_\varepsilon-\eta_\varepsilon'}{\varepsilon}b(t_\varepsilon,\eta_\varepsilon)-\frac{1}{2}\sigma^2(\eta_\varepsilon)A_{\varepsilon}\Big\}
           \leq 0\label{D},\\
     \end{split}\\
     \begin{split}
           \mbox{min}\Big\{ v^\Lambda(t_\varepsilon',\eta_\varepsilon')-\widehat{V}(t_\varepsilon',x_0-f(t_\varepsilon'),z_0,\eta_\varepsilon'),\frac{t_\varepsilon-t_\varepsilon'}{\varepsilon}
           -\frac{\eta_\varepsilon-\eta_\varepsilon'}{\varepsilon}b(t_\varepsilon',\eta_\varepsilon')-\frac{1}{2}\sigma^2(\eta_\varepsilon')B_{\varepsilon}\Big\}
           \geq  \Lambda\label{E},
     \end{split}
     \end{align}
where $b(t_\varepsilon,\eta_\varepsilon)=-\lambda(\eta_\varepsilon-\bar{\mu})$, $\sigma^2(\eta_\varepsilon)=\sigma_\mu^2$, $b(t_\varepsilon',\eta_\varepsilon')=-\lambda(\eta_\varepsilon'-\bar{\mu})$ and $\sigma^2(\eta_\varepsilon')=\sigma_\mu^2$.

    If $u-\widehat{V}\leq 0$ in (\ref{D}), then because $v^\Lambda-\widehat{V}\geq \Lambda$ in (\ref{E}), we obtain that $u-v^\Lambda\leq- \Lambda <0$ by contradiction with $\sup(u-v^\Lambda)=M>0$. On the other hand, if $u-\widehat{V}>0$ in (\ref{D}), then we have
    \begin{align*}
    \left\{
    \begin{array}{lr}
     \frac{t_\varepsilon-t_\varepsilon'}{\varepsilon}-\frac{\eta_\varepsilon-\eta_\varepsilon'}{\varepsilon}b(t_\varepsilon,\eta_\varepsilon)-\frac{1}{2}\sigma^2(\eta_\varepsilon)A_{\varepsilon}\leq 0,\\
     \frac{t_\varepsilon-t_\varepsilon'}{\varepsilon}-\frac{\eta_\varepsilon-\eta_\varepsilon'}{\varepsilon}b(t_\varepsilon',\eta_\varepsilon')-\frac{1}{2}\sigma^2(\eta_\varepsilon')B_{\varepsilon}
     \geq \Lambda.
    \end{array}
    \right.
    \end{align*}

    Furthermore, combining two inequalities above, we derive that
    \begin{align*}\begin{split}
     &\frac{\eta_\varepsilon-\eta_\varepsilon'}{\varepsilon}(b(t_\varepsilon,\eta_\varepsilon)-b(t_\varepsilon',\eta_\varepsilon'))
     +\frac{3}{2\varepsilon}(\sigma(\eta_\varepsilon)-\sigma(\eta_\varepsilon'))^2\\
     \geq
     &\frac{\eta_\varepsilon-\eta_\varepsilon'}{\varepsilon}(b(t_\varepsilon,\eta_\varepsilon)-b(t_\varepsilon',\eta_\varepsilon'))
     +\frac{1}{2}(\sigma^2(\eta_\varepsilon)A_{\varepsilon}-\sigma^2(\eta_\varepsilon')B_{\varepsilon})  \geq \Lambda.
    \end{split}\end{align*}
  The first inequality holds by Crandall-Ishii's lemma. In addition, by letting $\varepsilon\rightarrow 0$, we get $\frac{\eta_\varepsilon-\eta_\varepsilon'}{\varepsilon}(b(t_\varepsilon,\eta_\varepsilon)-b(t_\varepsilon',\eta_\varepsilon'))
     +\frac{3}{2\varepsilon}(\sigma(\eta_\varepsilon)-\sigma(\eta_\varepsilon'))^2=0$ thanks to (\ref{M}). It follows that we have $0\geq \Lambda>0$, which leads to a contradiction and therefore our claim holds.
\end{proof}

\begin{lemma}
For all $(t,\eta)\in\mathcal{C}$ in the continuation region, $\widetilde{V}$ in (\ref{tilde2v}) has H$\ddot{o}$lder continuous derivatives.
\end{lemma}

\begin{proof} The proof follows closely the argument in Section 6.3 of \cite{friedman2012stochastic}. First, let us recall that
 \begin{equation}\label{c2star3}
  \frac{\partial \widetilde{V}}{\partial t}(t,\eta) + \lambda(\eta-\bar{\mu})\frac{ \partial \widetilde{V} }{ \partial \eta }(t,\eta)
  -\frac{1}{2} \sigma_\mu^2 \frac{ \partial^2 \widetilde{V} }{ \partial \eta^2 }(t,\eta)=0 \ \mbox{on} \ \mathcal{C}.
 \end{equation}
The definition of viscosity solution of $\widetilde{V}$ to (\ref{varine2}) gives that $\widetilde{V}$ is a supersolution to (\ref{c2star3}). On the other hand, for any $(\bar{t},\bar{\eta})\in\mathcal{C}$, let $\varphi$ be a $C^2$ test function such that $(\bar{t},\bar{\eta})$ is a maximum of $\widetilde{V}-\varphi$ with $\widetilde{V}(\bar{t},\bar{\eta})=\varphi(\bar{t},\bar{\eta})$. By definition of $\mathcal{C}$, we have $\widetilde{V}(\bar{t},\bar{\eta}) > \widehat{V}(\bar{t},x_0-f(\bar{t}),z_0,\bar{\eta})$, so that
 \begin{equation*}
  \frac{\partial \varphi}{\partial t}(\bar{t},\bar{\eta}) + \lambda(\eta-\bar{\mu})\frac{ \partial \varphi }{ \partial \eta }(\bar{t},\bar{\eta})
  -\frac{1}{2} \sigma_\mu^2 \frac{ \partial^2 \varphi }{ \partial \eta^2 }(\bar{t},\bar{\eta})\leq0,
 \end{equation*}
due to the viscosity sub-solution property of $\widetilde{V}$ to (\ref{varine2}). It follows that $\widetilde{V}$ is a viscosity subsolution and therefore viscosity solution to (\ref{c2star3}).

Let us consider an initial boundary value problem:
 \begin{equation}\label{c2star2}
  \begin{split}
  -\frac{\partial w}{\partial t}(t,\eta) - \lambda(\eta-\bar{\mu})\frac{ \partial w }{ \partial \eta }(t,\eta)
  &+\frac{1}{2} \sigma_\mu^2 \frac{ \partial^2 w }{ \partial \eta^2 }(t,\eta)=0 \ \mbox{on} \ Q\cup B_T,    \\
  w(0,\eta)&=0 \ \mbox{on} \ B,\\
  w(t,\eta)&=\widehat{V}(t,x_0-\kappa t,z_0,\eta) \ \mbox{on} \ S.
  \end{split}
 \end{equation}
Here, $Q$ is an arbitrary bounded open region in $\mathcal{C}$, $Q$ lies in the strip $0<t<T$. $\tilde{B}=\bar{Q}\cap\{t=0\}$, $\tilde{B}_T=\bar{Q}\cap\{t=T\}$, $B_T$ denotes the interior of $\tilde{B}_T$, $B$ denotes the interior of $\tilde{B}$, $S_0$ denotes the boundary of $Q$ lying in the strip $0\leq t\leq T$ and $S=S_0\backslash B_T$. Theorem 3.6 in \cite{friedman2012stochastic} gives the existence and uniqueness of a solution $w$ on $Q\cup B_T$ to (\ref{c2star2}), and the solution $w$ has H$\ddot{o}$lder continuous derivatives $w_t$, $w_\eta$ and $w_{\eta\eta}$. Because the solution $w$ is a viscosity solution to (\ref{c2star3}) on $Q\cup B_T$, from standard uniqueness results on viscosity solution, we know that $\widetilde{V}=w$ on $Q\cup B_T$. As $Q\subset\mathcal{C}$ is arbitrary, it follows that $\widetilde{V}$ has the same property in the continuation region $\mathcal{C}$. Therefore, $\widetilde{V}$ has H\"{o}lder continuous derivatives $\widetilde{V}_t$, $\widetilde{V}_\eta$ and $\widetilde{V}_{\eta\eta}$.
\end{proof}
\ \\
Finally, we can prove Theorem \ref{unique}.
\begin{proof} We have shown the inequality $v^-=\sup_{p\in\mathcal{V}^-}p\leq \widetilde{V} \leq v^+=\inf_{q\in\mathcal{V}^+}q$ in Lemma \ref{compalem}. By using the comparison result in Proposition \ref{compp}, we also have $v^+\leq v^-$. Putting all pieces together, we conclude that $v^+=\widetilde{V}(t,\eta)=v^-$ and therefore the value function $\widetilde{V}(t,\eta)$ is the unique viscosity solution of the HJBVI \eqref{varine}. By following similar argument for Theorem 1 in \cite{Mihail}, fix the $\mathcal{F}_t$-adapted stopping time $\tau^*$ defined in \eqref{optimalentry}, It\^{o}-Tanaka's formula (see Theorem IV.1.5, Corollary IV.1.6 of \cite{Revvv}) can be applied to $\widetilde{V}(t,\mu_t)$ in view of H\"{o}lder continuous derivatives of $\widetilde{V}(t,\eta)$ and we get that
\begin{align*}
&\widehat{V}(\tau^*\wedge\tau_n, x_0-\kappa \tau^*\wedge\tau_n, z_0, \mu_{\tau^*\wedge\tau_n})\\
=&\widetilde{V}(t,\mu_t)+\left[ \widehat{V}(\tau^*\wedge\tau_n, x_0-\kappa \tau^*\wedge\tau_n, z_0, \mu_{\tau^*\wedge\tau_n})-\widetilde{V}(\tau^*\wedge\tau_n, \mu_{\tau^*\wedge\tau_n})\right]\\
+&\int_t^{\tau^*\wedge\tau_n}\sigma_{\mu}\frac{\partial \widetilde{V}}{\partial \eta}(s,\mu_s)dB_s+\int_t^{\tau^*\wedge\tau_n}\left[\frac{\partial \widetilde{V}(s,\mu_s)}{\partial t}+\mathcal{L}\widetilde{V}(s,\mu_s)\right]ds,
\end{align*}
where $\tau_n\uparrow T$ is the localizing sequence. As $\widetilde{V}(t,\eta)$ satisfies HJBVI \eqref{varine}, by taking conditional expectations and the definition of $\tau^*$ in \eqref{optimalentry}, we obtain that
\begin{align*}
\mathbb{E}_t\left[\widehat{V}(\tau^*\wedge\tau_n, x_0-\kappa \tau^*\wedge\tau_n, z_0, \mu_{\tau^*\wedge\tau_n})\mathbf{1}_{\{\tau^*\leq \tau_n\}}\right]+\mathbb{E}_t\left[\widetilde{V}(\tau_n, \mu_{\tau_n})\mathbf{1}_{\{\tau^*>\tau_n\}} \right] =\widetilde{V}(t,\mu_t)
\end{align*}
By taking the limit of $\tau_n$ and dominated convergence theorem, we can verify that
\begin{align*}
\mathbb{E}_t\left[\widehat{V}(\tau^*, x_0-\kappa \tau^*, z_0, \mu_{\tau^*})\right]=\widetilde{V}(t,\mu_t)
\end{align*}
and thus $\tau^*$ is the optimal entry time.

At last, the martingale property between $t=0$ and $\tau^*$ follows from the definition of stochastic subsolution and stochastic supersolution.
\end{proof}

Moreover, we can also easily verify the following sensitivity results of the composite value function.
\begin{lemma}
We have the following sensitivity properties of the value function $\widetilde{V}(t,\eta)$:

\begin{itemize}
\item[(i)] Suppose that $\alpha>$ and $\delta>0$ are both constants in the definition of habit formation process such that $\delta>\alpha$. We have that $\widetilde{V}(t,\eta; \alpha, \delta)$ is decreasing in $\delta$ and increasing in $\alpha$.
\item[(ii)] If the initial habit $z_0$ increases, the value function $\widetilde{V}(t,\eta)$ decreases.
\item[(iii)] If the information cost rate $\kappa$ increases, the value function $\widetilde{V}(t,\eta)$ decreases for any $t<T$.
\end{itemize}
\end{lemma}
\begin{proof}
By the definition of $\widetilde{V}(t,\eta)$ and the explicit form of $\widehat{V}(t, x_0-\kappa t, z_0, \eta)$ in \eqref{vtxz} and $m(t)$ in \eqref{mt}, for given $\delta>\alpha$, it is clear that $\widehat{V}(t, x_0-\kappa t, z_0, \eta)$ is decreasing in $\delta$ and increasing in $\alpha$, which implies that $\widetilde{V}(t,\eta)$ has the same sensitivity property. Similarly, it is clear that $\widehat{V}(t, x_0-\kappa t, z_0, \eta)$ decreases while $z_0$ increases, and hence $\widetilde{V}(t,\eta)$ is decreasing in $z_0$. At last, $\widehat{V}(t, x_0-\kappa t, z_0, \eta)$ decreases if $x_0-\kappa t$ decrease, it readily follows that $\widetilde{V}(t,\eta)$ is decreasing in $\kappa$.
\end{proof}

\appendix

\section{Explicit solution to the auxiliary ODEs}\label{appendix5.1}
Our ODE problems (\ref{atts}), (\ref{btts}), (\ref{ctts}) are similar to ODEs in \cite{brendle2006portfolio} for terminal wealth optimization problem, in which some insightful observations are made that we can solve these ODEs with time $t$ dependent coefficients by solving the following five auxiliary ODEs with constant coefficients, see Section 4 of \cite{brendle2006portfolio} for detailed discussions.

\begin{lemma}\label{lemma} For $k\leq t\leq s\leq T$, let us consider the following auxiliary ODEs for $a(t,s)$, $b(t,s)$, $l(t,s)$, $w(t,s)$ and $g(t,s)$:
    \begin{align}
      a_t=&-\frac{2(1-p+p\rho^2)}{1-p}\sigma^2_\mu a^2
      +\left(2\lambda-\frac{2p\rho\sigma_\mu}{(1-p)\sigma_S}\right)a
      -\frac{p}{2(1-p)\sigma^2_S}\label{lmaat},\\
      b_t=&-\frac{2(1-p+p\rho^2)}{1-p}\sigma^2_\mu ab
      -2\lambda\bar{\mu}a+\left(\lambda-\frac{p\rho\sigma_\mu}{(1-p)\sigma_S}\right)b\label{lmabt},\\
      l_t=&-\sigma^2_\mu a-\frac{(1-p+p\rho^2)\sigma^2_\mu}{2(1-p)}b^2-\lambda\bar{\mu}b\label{lmact},\\
      w_t=&-2(1-\rho^2)\sigma^2_\mu w^2+2\frac{\lambda\sigma_S+\rho\sigma_\mu}{\sigma_S}w+\frac{1}{2\sigma^2_S}\label{lmaft},\\
      g_t=&\sigma_\mu^2(1-\rho^2)(w-a)\label{lmagt}.
    \end{align}
    with the terminal conditions $a(s,s)=b(s,s)=l(s,s)=w(s,s)=g(s,s)=0$. Direct substitutions and computations show that the solutions of ODEs (\ref{atts}), (\ref{btts}), (\ref{ctts}) are given respectively by
    \begin{align}
      A(t,s)&:=\frac{a(t,s)}{(1-p)\left(1-2a(t,s)\hat{\Sigma}(t)\right)},\ \ \ \ B(t,s):=\frac{b(t,s)}{(1-p)\left(1-2a(t,s)\hat{\Sigma}(t)\right)},\notag\\
    \begin{split}
      C(t,s)&:=\frac{1}{1-p}\Big[l(t,s)+\frac{\hat{\Sigma}(t)}{\left(1-2a(t,s)\hat{\Sigma}(t)\right)}b^2(t,s)
      -\frac{1-p}{2}\log\left(1-2a(t,s)\hat{\Sigma}(t)\right)\\
      &-\frac{p}{2}\log\left(1-2w(t,s)\hat{\Sigma}(t)\right)-pg(t,s)\Big].\label{preODEs}
    \end{split}
    \end{align}
\end{lemma}

Following the same arguments in \cite{kim1996dynamic}, we can actually solve auxiliary ODEs (\ref{lmaat}), (\ref{lmabt}), (\ref{lmact}), (\ref{lmaft}) and (\ref{lmagt}) explicitly in the order that we first solve the simple ODEs (\ref{lmaat}) and (\ref{lmaft}) to get $a(t,s)$ and $w(t,s)$, and then obtain $b(t,s)$ and $g(t,s)$ by solving ODEs (\ref{lmabt}) and (\ref{lmagt}). At last, we solve ODE (\ref{lmact}) to conclude $l(t,s)$. We therefore can get that

\begin{align*}
      a(t,s)=&\frac{p(1-e^{2\xi(t-s)})}
      {2(1-p)\sigma_S^2\Big[2\xi-(\xi+\gamma_2)(1-e^{2\xi(t-s)})\Big]},\\
      b(t,s)=&\frac{p\lambda\bar{\mu}(1-e^{\xi(t-s)})^2}
      {(1-p)\sigma_S^2\xi\Big[2\xi-(\xi+\gamma_2)(1-e^{2\xi(t-s)})\Big]},\\
    \begin{split}
        l(t,s)=&\frac{p}{2(1-p)\sigma_S^2}\left(\frac{\lambda^2\bar{\mu}^2}{\xi^2}-\frac{\sigma_\mu^2\gamma_2}{\gamma_2^2-\xi^2}\right)(s-t)\\
        &+\frac{p\lambda^2\bar{\mu}^2\Big[(\xi+2\gamma_2)e^{2\xi(t-s)}-4\gamma_2e^{\xi(t-s)}+2\gamma_2-\xi\Big]}
        {2(1-p)\sigma_S^2\xi^3\left[2\xi-(\xi+\gamma_2)(1-e^{2\xi(t-s)})\right]}\\
        &+\frac{p\sigma_\mu^2}{2(1-p)\sigma_S^2(\xi^2-\gamma_2^2)}\log\left|\frac{2\xi-(\xi+\gamma_2)(1-e^{2\xi(t-s)})}{{2\xi}e^{\xi(t-s)}}\right|,
    \end{split}\\
      w(t,s)=&-\frac{1}{2\sigma_S}\frac{1-e^{2\xi_1(t-s)}}
      {(\sigma_S\xi_1+\lambda\sigma_S+\rho\sigma_\mu)+(\sigma_S\xi_1-\lambda\sigma_S-\rho\sigma_\mu)e^{2\xi_1(t-s)}},\\
    \begin{split}
        g(t,s)=&\frac{1}{2}\log\left(\frac{(\sigma_S\xi_1+\lambda\sigma_S+\rho\sigma_\mu)+(\sigma_S\xi_1-\lambda\sigma_S-\rho\sigma_\mu)e^{2\xi_1(t-s)}}{2\sigma_S\xi_1e^{\xi_1(t-s)}}\right)\\
        &-\frac{(1-p)(1-\rho^2)}{2(1-p+p\rho^2)}\log\left(\frac{(\sigma_S\xi+\lambda\sigma_S-\frac{\rho\sigma_\mu p}{1-p})+(\sigma_S\xi-\lambda\sigma_S+\frac{\rho\sigma_\mu p}{1-p})e^{2\xi(t-s)}}{2\sigma_S\xi e^{\xi(t-s)}}\right)\\
        &-\frac{\rho^2\lambda(s-t)}{2(1-p+p\rho^2)}-\frac{\rho\sigma_\mu(s-t)}{2(1-p+p\rho^2)\sigma_S},
    \end{split}
    \end{align*}
    where
    \begin{align}\label{delta}
      \Delta:=\lambda^2-\frac{2\lambda p\rho\sigma_\mu}{(1-p)\sigma_S}-\frac{p\sigma_\mu^2}{(1-p)\sigma_S^2}>0,
    \end{align}
and
    \begin{align*}
      \xi:=\sqrt\Delta=\sqrt{\gamma^2_2-\gamma_1\gamma_3},\ \
      \xi_1:=\frac{\sqrt{(1-\rho^2)\sigma_\mu^2+(\lambda\sigma_S+\rho\sigma_\mu)^2}}{\sigma_S},\\
      \gamma_1:=\frac{(1-p+p\rho^2)}{1-p}\sigma_\mu^2, \ \ \gamma_2:=-\lambda+\frac{p\rho\sigma_\mu}{(1-p)\sigma_S}, \ \ \gamma_3:=\frac{p}{(1-p)\sigma_S^2}.
    \end{align*}

Moreover, it is straightforward to see that $a$, $b$, $l$, $w$ and $g$ are globally bounded if we have that $\gamma_3>0$, or  $\gamma_1>0$, or $\gamma_2<0$.

\begin{remark}\label{appendixrk1}
Under the assumption that $p<0$, (\ref{delta}) clearly holds and we have $\gamma_2<0$. We can see that $a(t,s)\leq 0$ and $b(t,s)\leq 0$ are bounded and $1-2a(t,s)\hat{\Sigma}(t)>1$ and $1-w(t,s)\hat{\Sigma}(t)>1$. By expressions in \eqref{preODEs}, we can conclude that $A(t,s)$, $B(t,s)$ and $C(t,s)$ are all bounded on $k\leq t\leq s\leq T$ as well as $A(t,s)=\frac{a(t,s)}{(1-p)(1-2a(t,s)\hat{\Sigma}(t))}\leq 0$ and $B(t,s)=\frac{b(t,s)}{(1-p)\left(1-2a(t,s)\hat{\Sigma}(t)\right)}\leq 0$, for $k\leq t\leq s\leq T$.
\end{remark}

\section{Proof of the verification theorem}\label{appendix5.2}
We first show that the consumption constraint $c_t\geq Z_t$ implies the constraint on the controlled wealth process in the next lemma.
\begin{lemma}\label{lma4pre}
The admissible space $\mathcal{A}$ is not empty if and only if the initial budget constraint
$x\geq m(k)z$ is fulfilled. Moreover, for each pair $(\pi,c)\in\mathcal{A}$, the controlled wealth process $\hat{X}^{\pi,c}_t$ satisfies the constraint
\begin{align}\label{xpict}
  \hat{X}^{\pi,c}_t\geq m(t)Z_t, \ \ k\leq t\leq T,
\end{align}
where the deterministic function $m(t)$ is defined in (\ref{mt}) and refers to the cost of subsistence consumption per unit of standard of living at time t.
\end{lemma}
\begin{proof}
Let's first assume that $x\geq m(k)z$, we can always take $\pi_t\equiv 0$, and $c_t=ze^{\int^t_k(\delta(v)-\alpha(v))dv}$ for $t\in[k,T]$. It is easy to verify $\hat{X}^{\pi,c}_t\geq 0$ and $c_t\equiv Z_t$ so that $(\pi,c)\in\mathcal{A}$, and hence $\mathcal{A}$ is not empty.

    On the other hand, starting from $t=k$ with the wealth $x$ and the standard of living $z$, the addictive habits constraint $c_t\geq Z_t$, $k\leq t\leq T$ implies that the consumption must always exceed the \emph{subsistence consumption} $\bar{c}_t=Z(t;\bar{c}_t)$ which satisfies
    \begin{align}\label{dbarct}
      d\bar{c}_t=(\delta(t)-\alpha(t))\bar{c}_tdt, \ \ \bar{c}_k=z, \ \ k\leq t\leq T.
    \end{align}
    Indeed, since $Z_t$ satisfies $dZ_t=(\delta_tc_t-\alpha_tZ_t)dt$ with $Z_k=z\geq 0$, the constraint $c_t\geq Z_t$ implies that
    \begin{align}\label{dzt}
      dZ_t\geq(\delta_tZ_t-\alpha_tZ_t)dt, \ \ Z_k=z.
    \end{align}
    By (\ref{dbarct}) and (\ref{dzt}), one can get $d(Z_t-\bar{c}_t)\geq(\delta_t-\alpha_t)(Z_t-\bar{c}_t)dt$ and $Z_k-\bar{c}_k=0$, from which we can derive that $e^{\int^t_k(\delta_s-\alpha_s)ds}(Z_t-\bar{c}_t)\geq 0$, $k\leq t\leq T$. It follows that $c_t\geq\bar{c}_t$, which is equivalent to
    \begin{align}\label{ct}
      c_t\geq ze^{\int^t_k(\delta(v)-\alpha(v))dv}, \ \ k\leq t\leq T.
    \end{align}

    Define the exponential local martingale $\widetilde{H}_t=\exp\left(-\int^t_k\frac{\hat{\mu}_v}{\sigma_S}d\hat{W}_v-\frac{1}{2}\int^t_k\frac{\hat{\mu}_v^2}{\sigma_S^2}dv\right)$, $k\leq t\leq T$.
    As $\hat{\mu}_t$ follows the dynamics (\ref{dhatmut}), we derive that
    \begin{align*}
      \hat{\mu}_t=e^{-t\lambda}\eta+\bar{\mu}(1-e^{-t\lambda})
      +\int^t_ke^{\lambda(u-t)}\frac{\left(\hat{\Sigma}(u)+\sigma_S\sigma_\mu\rho\right)}{\sigma_S}d\hat{W}_u.
    \end{align*}
    Similar to the proof of Corollary 3.5.14 and Corollary 3.5.16 in \cite{motion1988brownian}, Bene\v s' condition implies that $\widetilde{H}$ is a true martingale with respect to $(\Omega,\mathcal{F}^S,\mathbb{P})$.

    Now, define the probability measure $\widetilde{\mathbb{P}}$ as $ \frac{d\widetilde{\mathbb{P}}}{d\mathbb{P}}=\widetilde{H}_T$, Girsanov theorem states that $\widetilde{W}_t:=\hat{W}_t+\int^t_k\frac{\hat{\mu}_v}{\sigma_S}dv$, $k\leq t\leq T$ is a Brownian Motion under $(\widetilde{\mathbb{P}},(\mathcal{F}^S_t)_{k\leq t\leq T})$. We can rewrite the wealth process as $\hat{X}_T+\int^T_kc_vdv=x+\int^T_k\pi_v\sigma_Sd\widetilde{W}_v$.
    As we have $\hat{X}_T\geq 0$, it is easy to see that $\int^t_k\pi_v\sigma_Sd\widetilde{W}_v$ is a supermartingale under $(\Omega,\mathbb{F}^S,\widetilde{\mathbb{P}})$. By taking the expectation under $\widetilde{\mathbb{P}}$, we have $x\geq\widetilde{\mathbb{E}}\left[\int^T_kc_vdv\right]$. Thanks to the inequality (\ref{ct}), we further have $
      x\geq z\widetilde{\mathbb{E}}\left[\int^T_k\exp\left(\int^v_k(\delta(u)-\alpha(u))du\right)dv\right]$.
    Because $\delta(t)$ and $\alpha(t)$ are deterministic functions, we obtain that $x\geq m(k)z$. In general, for $\forall t\in[k,T]$, following the same procedure, we can take conditional expectation under filtration $\mathcal{F}^S_t$, and get
      $\hat{X}_t\geq Z_t\widetilde{\mathbb{E}}\bigg[\int^T_t\exp\Big(\int^v_t(\delta(u)-\alpha(u))du\Big)dv\bigg{|}\mathcal{F}^S_t\bigg]$.
    Again as $\delta(t)$, $\alpha(t)$ are deterministic, we get $\hat{X}_t\geq m(t)Z_t$, $k\leq t\leq T$.
\end{proof}
\ \\
We can finally prove Theorem \ref{thm} for the interior control problem.
\begin{proof} For any pair of admissible control $(\pi_t,c_t)\in\mathcal{A}$, It\^o's lemma gives
\begin{align}\label{dV}
  d\left[V(t,\hat{X}_t,Z_t,\hat{\mu}_t)\right]=\left[\mathcal{G}^{\pi_t,c_t}V(t,\hat{X}_t,Z_t,\hat{\mu}_t)\right]dt
  +\left[V_{x}\sigma_S\pi_t+V_{\eta}\frac{\left(\hat{\Sigma}(t)+\sigma_S\sigma_\mu\rho\right)}{\sigma_S}\right]d\hat{W}_t,
\end{align}
where we define the process $\mathcal{G}^{\pi_t,c_t}V(t,\hat{X}_t,Z_t,\hat{\mu}_t)$ by
\begin{align}
\begin{split}
&\mathcal{G}^{\pi_t,c_t}V(t,\hat{X}_t,Z_t,\hat{\mu}_t)=V_t-\alpha(t)Z_tV_t-\lambda(\hat{\mu}_t-\bar{\mu})V_{\eta}
  +\frac{\left(\hat{\Sigma}(t)+\sigma_S\sigma_\mu\rho\right)^2}{2\sigma_S^2}V_{\eta\eta}-c_tV_{x}\\
 &+c_t\delta(t)V_{z}
  +\frac{(c_t-Z_t)^p}{p}+\pi_t\hat{\mu}_tV_{x}+\frac{1}{2}\sigma_S^2\pi_t^2V_{xx}
  +V_{x\eta}\left(\hat{\Sigma}(t)+\sigma_S\sigma_\mu\rho\right)\pi_t.\nonumber
\end{split}
\end{align}
For any localizing sequence $\tau_n$, by integrating (\ref{dV}) on $[k,\tau_n\wedge T]$ and taking the expectation, we have
    \begin{align}\label{v0x0z0eta0}
      V(k,x,z,\eta)\geq\mathbb{E}\left[\int^{\tau_n\wedge T}_k\frac{(c_s-Z_s)^p}{p}ds\right]
        +\mathbb{E}\Big[V(\tau_n\wedge T,\hat{X}_{\tau_n\wedge T},Z_{\tau_n\wedge T},\hat{\mu}_{\tau_n\wedge T})\Big].
    \end{align}

    Similar to the argument in \cite{janevcek2012optimal}, let us consider a fixed pair of control $(\pi_t,c_t)\in\mathcal{A}=\mathcal{A}_{x}$, where we denote $\mathcal{A}_{x}$ as the admissible space with initial endowment $x$. For $\forall\epsilon>0$, it is clear that $\mathcal{A}_{x}\subseteq\mathcal{A}_{x+\epsilon}$, and $(\pi_t,c_t)\in\mathcal{A}_{x+\epsilon}$. Also it is easy to see that $\hat{X}^{x+\epsilon}_t=\hat{X}^{x}_t+\epsilon=\hat{X}_t+\epsilon$, $k\leq t\leq T$. As the process $Z_t$ is defined using this consumption policy $c_t$, under the probability measure $\mathbb{P}_{x,z,\eta}$, we can obtain
    \begin{align}\label{v0x0eps}
      V(k,x+\epsilon,z,\eta)\geq\mathbb{E}\left[\int^{\tau_n\wedge T}_k\frac{(c_s-Z_s)^p}{p}ds\right]
        +\mathbb{E}\Big[V(\tau_n\wedge T,\hat{X}_{\tau_n\wedge T}+\epsilon,Z_{\tau_n\wedge T},\hat{\mu}_{\tau_n\wedge T})\Big].
    \end{align}

    Monotone Convergence Theorem first leads to
    \begin{align*}
      \lim_{n\to+\infty}\mathbb{E}\left[\int^{\tau_n\wedge T}_k\frac{(c_s-Z_s)^p}{p}ds\right]
        =\mathbb{E}\left[\int^{T}_k\frac{(c_s-Z_s)^p}{p}ds\right].
    \end{align*}

    For simplicity, let's denote $Y_t=\Big(\hat{X}_t-m(t)Z_t\Big)$. The definition (\ref{vtxz}) implies that: $
      V(\tau_n\wedge T,\hat{X}_{\tau_n\wedge T}+\epsilon,Z_{\tau_n\wedge T},\hat{\mu}_{\tau_n\wedge T})=\frac{1}{p}(Y_{\tau_n\wedge T}+\epsilon)^pN^{1-p}_{\tau_n\wedge T}$.
    Lemma \ref{lma4pre} gives $\hat{X}_t\geq m(t)Z_t$ for $k\leq t\leq T$ under any admissible control $(\pi_t,c_t)$, we get that $Y_{\tau_n\wedge T}+\epsilon\geq\epsilon>0$, $\forall k\leq t\leq T$. As $p<0$, it follows that
    \begin{align}\label{sup}
      \sup_{n}(Y_{\tau_n\wedge T}+\epsilon)^p<\epsilon^p<+\infty.
    \end{align}

    Remark \ref{appendixrk1} gives that $A(t,s)\leq 0$, $B(t,s)$ and $C(t,s)$ are all bounded on $k\leq t\leq s\leq T$, $\forall k\leq t\leq s\leq T$. Also $m(s)$, $\delta(s)$ are continuous functions and hence bounded on $[k,T]$. Hence $N(k,\eta)\leq k_1\exp(k_2\eta)$, for some constants $k_2, k_1>1$.
    It follows that there exist some constants $\bar{k}_2,\bar{k}_1>1$ such that
    \begin{align*}
      \sup_{n}N^{1-p}_{\tau_n\wedge T}\leq\sup_{t\in[k,T]}\Big(k_1\exp(k_2\hat{\mu}_t)\Big)^{1-p}
      \leq\bar{k}_1\exp\left(\bar{k}_2\sup_{t\in[k,T]}\hat{\mu}_t\right).
    \end{align*}

    The process $\hat{\mu}_t$ satisfies (\ref{dhatmut}), which leads to
    \begin{align*}
      \hat{\mu}_t=e^{-t\lambda}\eta+\bar{\mu}(1-e^{-t\lambda})
        +\int^t_ke^{\lambda(u-t)}\frac{\left(\hat{\Sigma}(u)+\sigma_S\sigma_\mu\rho\right)}{\sigma_S}d\hat{W}_u.
    \end{align*}

    Hence, there exists positive constants $l$ and $l_1>1$ large enough, such that
      $\sup_{t\in[k,T]}\hat{\mu}_t\leq l+\sup_{t\in[k,T]}l_1\hat{W}_t$, $t\in[k,T]$. Using the distribution of running maximum of the Brownian Motion, there exist some positive constants $\bar{l}>1$ and $\bar{l}_1$ such that
    \begin{align}\label{esup}
      \mathbb{E}\left[\sup_{n}N^{1-p}_{\tau_n\wedge T}\right]
        \leq\bar{l}_1\mathbb{E}\left[\exp(\sup_{t\in[k,T]}\bar{l}\hat{B}_t)\right]<+\infty.
    \end{align}

    At last, by (\ref{sup}) and (\ref{esup}), we can conclude that
    \begin{align*}
    \mathbb{E}\left[\sup_{n}V(\tau_n\wedge T,\hat{X}_{\tau_n\wedge T}+\epsilon,Z_{\tau_n\wedge T},\hat{\mu}_{\tau_n\wedge T})\right]<+\infty.
    \end{align*}
    Dominated Convergence Theorem and $N(T,\hat{\mu_T}) = 0$ imply that
    \begin{align*}
      \lim_{n\to\infty}\mathbb{E}\Big[V(\tau_n\wedge T,\hat{X}_{\tau_n\wedge T}+\epsilon,Z_{\tau_n\wedge T},\hat{\mu}_{\tau_n\wedge T})\Big]=\mathbb{E}\left[\frac{1}{p}(Y_T+\epsilon)^pN^{1-p}(T,\hat{\mu_T})\right] = 0.
    \end{align*}
    Combining this with equation (\ref{v0x0eps}) and $(\pi_t,c_t)\in\mathcal{A}$, we have that
    \begin{align*}
      V(k,x+\epsilon,z,\eta;\theta)\geq\sup_{\pi,c\in\mathcal{A}}\mathbb{E}
        \left[\int^T_k\frac{(c_s-Z_s)^p}{p}ds\right]
        =\widehat{V}(k,x,z,\eta,\theta).
    \end{align*}
    Note that $V(t,x,z,\eta;\theta)$ is continuous in variable $x$. By letting $\epsilon\to 0$, we deduce that
    \begin{align*}
      V(k,x,z,\eta;\theta)=\lim_{\epsilon\to 0}V(k,x+\epsilon,z,\eta)
        \geq \widehat{V}(k,x,z,\eta,\theta).
    \end{align*}

    On the other hand, for $\pi^*_t$ and $c^*_t$ given in (\ref{pistar}) and (\ref{cstar}), we first need to show that the SDE
    \begin{align}\label{dxstart}
      d\hat{X}^*_t=(\pi^*_t\mu_t-c^*_t)dt+\sigma_S\pi^*_td\hat{W}_t, \ \ k\leq t\leq T,
    \end{align}
    with initial condition $x>m(k)z$ admits a unique strong solution that satisfies the constraint $\hat{X}^*_t>m(t)Z^*_t$, $\forall k\leq t\leq T$. Denote $Y^*_t=\hat{X}^*_t-m(t)Z^*_t$. By It\^o's lemma and substitution of $c^*_t$ using (\ref{cstar}), we obtain that
    \begin{align*}
    \begin{split}
      dY^*_t=&\left[-\frac{\Big(1+\delta(t)m(t)\Big)^{\frac{-p}{1-p}}}{N}
             +\frac{\hat{\mu}_t^2}{(1-p)\sigma_S^2}+\frac{\Big(\hat{\Sigma}(t)+\sigma_S\sigma_\mu\rho\Big)}{\sigma_S^2}\frac{N_{\eta}}{N}\hat{\mu}_t\right]Y^*_t dt\\
             &+\left[\frac{\hat{\mu}_t}{(1-p)\sigma_S}+\frac{\Big(\hat{\Sigma}(t)+\sigma_S\sigma_\mu\rho\Big)}{\sigma_S}\frac{N_{\eta}}{N}\right]Y^*_td\hat{W}_t.
    \end{split}
    \end{align*}

    Let us define the auxiliary process $\Gamma_t := \frac{N(t,\hat{\mu}_t)}{Y^*_t}$, for $k\leq t\leq T$. It\^o's lemma gives that
    \begin{align}
    \begin{split}\label{dgammat}
      d\Gamma_t =&\frac{\Gamma_t}{N_t}\Bigg[N_t-\lambda(\hat{\mu}_t-\bar{\mu})N_{\eta}
                +\frac{\left(\hat{\Sigma}(t)+\sigma_S\sigma_\mu\rho\right)^2}{2\sigma_S^2}N_{\eta\eta}
                +\frac{\hat{\mu}_t\left(\hat{\Sigma}(t)+\sigma_S\sigma_\mu\rho\right)p}{(1-p)\sigma_S^2}N_{\eta}\\
                &+\Big(1+\delta(t)m(t)\Big)^{\frac{-p}{1-p}}
                +\frac{p\hat{\mu}_t^2}{(1-p)^2\sigma_S^2}N\Bigg]dt
                +\Gamma_t\left[\frac{-\hat{\mu}_t}{(1-p)\sigma_S}\right]d\hat{W}_t.
    \end{split}
    \end{align}

    As $N(t,\eta)$ satisfies the linear PDE (\ref{nteta}), (\ref{dgammat}) is reduced to $d\Gamma_t=\Gamma_t\left[\frac{p\hat{\mu}_t^2}{2(1-p)^2\sigma_S^2}\right]dt
               +\Gamma_t\left[\frac{-\hat{\mu}_t}{(1-p)\sigma_S}\right]d\hat{W}_t$, and the existence of the unique strong solution is verified and $\Gamma_k=\frac{N(k,\eta)}{x-m(k)z}>0$ implies that $\Gamma_t>0$, $\forall k\leq t\leq T$. It holds that the SDE (\ref{dxstart}) admits a unique strong solution defined in (\ref{xstart}) and the solution $\hat{X}^*_t$ satisfies the constraint (\ref{xpict}).

    Next, we verify the pair $(\pi^*_t,c^*_t)$ is indeed in the admissible space $\mathcal{A}$. First, by the definition (\ref{pistar}) and (\ref{cstar}), it is clear that $\pi^*_t$ and $c^*_t$ are $\mathcal{F}^S_t$ progressively measurable, and by the path continuity of $Y^*_t=\hat{X}^*_t-m(t)Z^*_t$ and of $\pi^*_t$ and $c^*_t$, it is easy to show that $\int^T_k(\pi^*_t)^2dt<+\infty$ and $\int^T_kc^*_tdt<+\infty$, a.s. Also, because $\hat{X}^*_t>m(t)Z^*_t$, $\forall t\in[k,T]$, by the definition of $c^*_t$, the consumption constraint $c^*_t>Z^*_t$, $\forall t\in[k,T]$ is satisfied. It follows that $(\pi^*_t,c^*_t)\in\mathcal{A}$. Given $(\pi^*_t,c^*_t)$ as above, the equality is proved that
    \begin{align*}
      V(k,x,z,\eta;\theta)=\mathbb{E}\left[\int^{\tau_n\wedge T}_k\frac{(c^*_t-Z^*_t)^p}{p}dt\right]
        +\mathbb{E}\Big[V(\tau_n\wedge T,\hat{X}^*_{\tau_n\wedge T},Z^*_{\tau_n\wedge T},\hat{\mu}_{\tau_n\wedge T})\Big].
    \end{align*}
    Monotone Convergence Theorem gives $\lim_{n\to+\infty}\mathbb{E}\left[\int^{\tau_n\wedge T}_k\frac{(c^*_t-Z^*_t)^p}{p}dt\right]
        =\mathbb{E}\left[\int^{T}_k\frac{(c^*_t-Z^*_t)^p}{p}dt\right]$. Moreover, as we have $V(t,x,z,\eta)<0$ by $p<0$, Fatou's lemma implies that
    \begin{align*}
      \limsup_{n\to+\infty}\mathbb{E}\Big[V(\tau_n\wedge T,\hat{X}^*_{\tau_n\wedge T},Z^*_{\tau_n\wedge T},\hat{\mu}_{\tau_n\wedge T})\Big]
        \leq\mathbb{E}\Big[V(T,\hat{X}^*_{T},Z^*_{T},\hat{\mu}_{T})\Big]
        =0.
    \end{align*}

    It follows that $V(k,x,z,\eta;\theta)
        \leq\mathbb{E}\left[\int^{T}_k\frac{(c^*_t-Z^*_t)^p}{p}dt\right]
        \leq \widehat{V}(k,x,z,\eta,\theta)$, which completes the proof.
\end{proof}


%
%
%
%


\ \\
\textbf{Acknowledgements} Y. Yang and X. Yu are supported by the Hong Kong Polytechnic University research grant under no. P0031417.

\end{document}